%% file: paper.tex
\documentclass{siamltex}

\usepackage{amsmath}
\usepackage{import}
\usepackage{amssymb}
\usepackage{amsfonts}
\usepackage{amsbsy}
\usepackage{fancyvrb}
\usepackage{stmaryrd}
\usepackage{pdfsync}
\usepackage{graphicx}
\usepackage{algorithm}
\usepackage{color}
\usepackage[shadow]{todonotes}
\usepackage[numbers,square,sort&compress]{natbib}

\DefineVerbatimEnvironment{code}{Verbatim}{frame=single,rulecolor=\color{blue}}
\newcommand{\tab}{\hspace*{2em}}
\newcommand{\codesize}{\footnotesize}
\newcommand{\mesh}{\mathcal{T}}

\newcommand{\emp}[1]{\texttt{#1}}
\newcommand{\OmegaD}{\Omega_{0,\mathrm{D}}}
\newcommand{\OmegaN}{\Omega_{0,\mathrm{N}}}
\newcommand{\avg}[1]{\langle #1 \rangle}

\input{mydefinitions.sty}

\title{Efficient implementation of finite element methods on
  non-matching and overlapping meshes in 3D}

\author{Andr\'e Massing \and Mats G. Larson \and Anders Logg}

\begin{document}

\maketitle

\begin{abstract}
  In recent years, a number of finite element methods have been
  formulated for the solution of partial differential equations on
  complex geometries based on non-matching or overlapping
  meshes. Examples of such methods include the fictitious domain
  method, the extended finite element method, and Nitsche's method. In
  all of these methods, integrals must be computed over cut cells or
  subsimplices which is challenging to implement, especially in three
  space dimensions. In this note, we address the main challenges of
  such an implementation and demonstrate good performance of a fully
  general code for automatic detection of mesh intersections and
  integration over cut cells and subsimplices. As a canonical example
  of an overlapping mesh method, we consider Nitsche's method which we
  apply to Poisson's equation and a linear elastic problem.
\end{abstract}


\section{Introduction}

A fundamental problem in computational science is to develop methods
for the solution of partial differential equations on domains
containing one or several objects that may have complex or
time-dependent geometry. One approach to attacking this problem is to
allow overlapping meshes where a mesh representing an object is
allowed to overlap a background mesh representing the surroundings of
the object; see for instance \citet{ChesshireHenshaw1990,Yu2005,
  ZhangGerstenbergerWangEtAl2004, MayerPoppGerstenbergerEtAl2010} for
various applications.

In solid mechanics, overlapping meshes may be used to represent
materials consisting of elastic objects inserted into a surrounding
elastic material of another type \citep{Jirasek2002}, and in
fluid--structure interaction, an overlapping mesh may be used to
represent an elastic body immersed in a fluid represented by a fixed
background mesh
\citep{Yu2005,BaigesCodina2009,MayerGerstenbergerWall2009}. Another
common application \citep{ChesshireHenshaw1990,Dhia2005,DayBochev2008}
is found in mesh generation where a complicated geometry such as, for
example, a pipe junction, may be decomposed into simpler parts and one
unstructured tetrahedral mesh is created for each part.  These
components may then be stored, reused and recombined in applications
by using an overlapping mesh technique.

Overlapping mesh techniques are of particular interest in simulations
that involve moving objects. For such problems, overlapping mesh
techniques are an attractive alternative to ALE techniques. The main
advantage is that by using an overlapping mesh technique, one avoids
deformation of the mesh that may lead to deterioration of the mesh
quality and ultimately force remeshing. This is of particular
importance in the simulation of fluid--structure interaction where the
topology of the fluid domain may change due to deformation of the
solid.

The main focus of this work is on the general algorithms and efficient
implementation that is required to handle complex problems posed on
overlapping meshes in three space dimensions. Many of the presented
algorithms and the tools developed are of interest and use for the
implementation of various overlapping mesh techniques. To make the
discussion concrete, we here focus mainly on Nitsche's method. In
~\citet{HansboHansboLarson2003}, a consistent finite element method
for overlapping meshes based on Nitsche's method was introduced and
analyzed. The basic idea is to construct a finite element space by
taking the direct sum of the space of continuous piecewise polynomial
functions on the overlapping mesh and the restriction of the space of
continuous piecewise polynomial functions to the complement of the
overlapping mesh, and then impose the interface conditions using
Nitsche's method. It was shown that this approach leads to a stable
method of optimal order for arbitrary degree polynomial approximation.

The main challenge in the implementation is to compute the
intersection between the overlapping and the background mesh. The
result is a set of \emph{cut} cells which may be arbitrarily complex
polyhedra. These arise as the result of subtracting from the
tetrahedra of the background mesh a set of tetrahedra from the
overlapping mesh. By adopting algorithms and search structures from
the field of computational geometry, we show how these issues can be
handled in an efficient manner. Furthermore, one needs to compute
integrals on the resulting polyhedra. This can be carried out
efficiently based on an application of the divergence theorem in
combination with potentials to represent an integral on the
three-dimensional polyhedron as a sum of one-dimensional integrals on
its edges.

The presented algorithms and implementation are relevant for several
other types of related methods, including the extended finite element
method (XFEM) \citep{SukumarMoesMoranEtAl2000,FriesBelytschko2010}, non-fitted
sharp interface methods
\citep{HansboHansbo2002,BeckerBurmanHansbo2009}, and mesh-tying
techniques \citep{DayBochev2008}.

\subsection{Major contributions of this paper}

Our work consists of several contributions. We identify the major
techniques used in the implementation of overlapping mesh methods and
related methods. We further identify useful data structures and
algorithms from the field of computational geometry. As part of our
work, existing computational geometry libraries such as
\texttt{CGAL}~\cite{cgal} and \texttt{GTS}~\cite{gts} have been
wrapped into the general purpose finite element library
DOLFIN~\cite{LoggWells2010a} and into an extension library on top of
it, thereby making these algorithms and data structures more easily
accessible to the finite element community. Based on our
implementation, we demonstrate for the first time a highly efficient
implementation of Nitsche's method on overlapping meshes for several
problems posed in three spatial dimensions, thus opening the
possibility of employing Nitsche-based overlapping mesh methods for
challenging 3D problems such as fluid--structure interaction or
domain-bridging problems.

\subsection{Outline of this paper}

In Section 2, we review Nitsche's overlapping mesh method for a model
problem and present in Section~3 the techniques and algorithms we have
developed for efficient implementation of Nitsche's method in three
space dimensions. The corresponding implementation and data structures
are described in Section~4. Finally, we present in Section~5 numerical
examples that demonstrate the convergence of the numerical solution as
well as the scaling of the work required to compute mesh intersections
relative to standard finite element assembly on matching meshes.

\section{Nitsche's method on overlapping meshes}

We here review Nitsche's method for a simple model problem posed on
two overlapping meshes.

\subsection{Model problem}

\begin{figure}
  \begin{center}
    \includegraphics[width=0.40\textwidth]{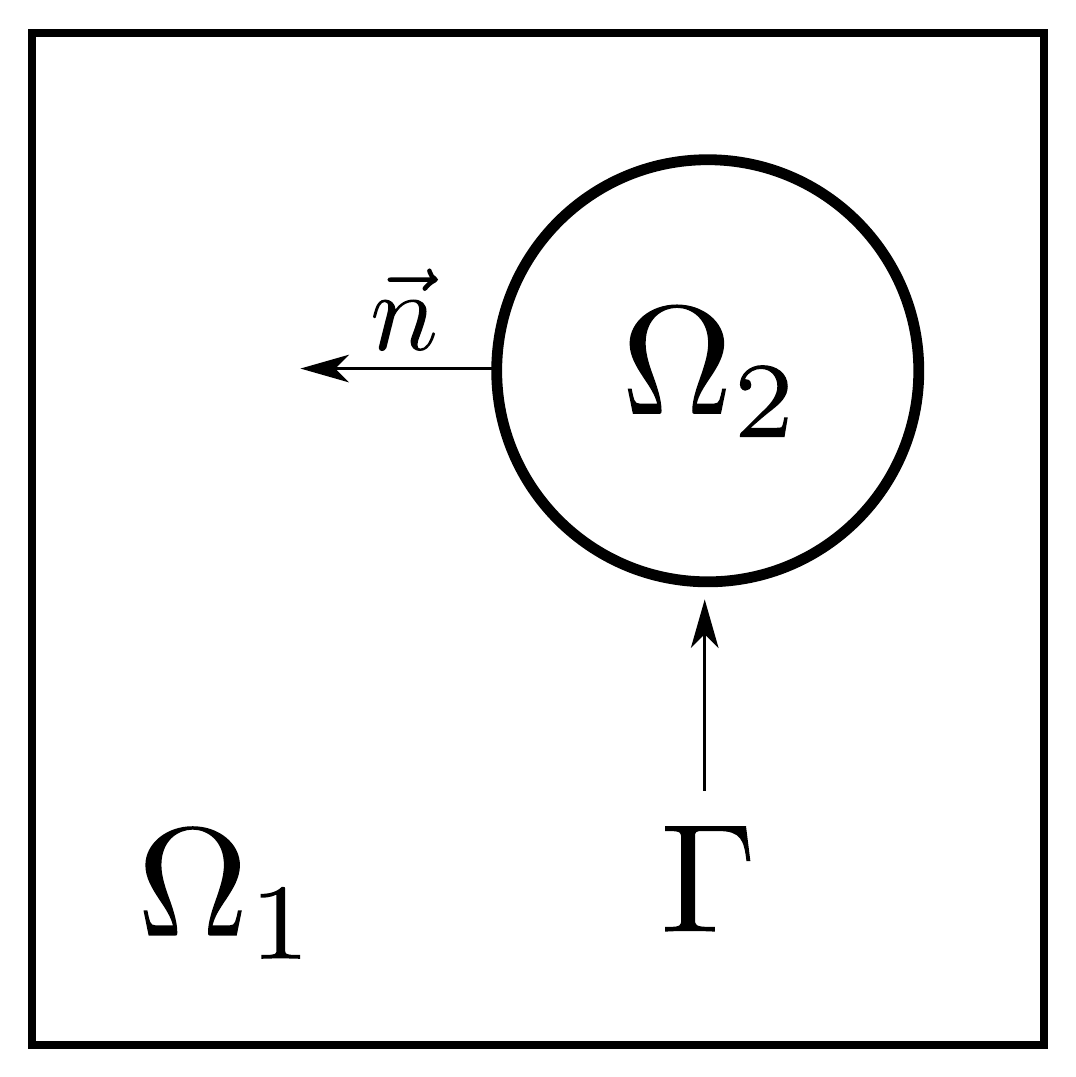}
    \caption{Two domains $\Omega_1$ and $\Omega_2$ separated by the
      common interface $\Gamma = \partial \Omega_1 \cap \partial
      \Omega_2$.}
  \label{fig:model_domain}
  \end{center}
\end{figure}

Let $\Omega_0 = (\overline{\Omega_1 \cup \Omega_2})^{\circ}$ be a
domain in $\R^3$, with boundary $\partial \Omega_0$, consisting of two
(open and bounded) subdomains $\Omega_1$ and $\Omega_2$ separated by
the interface $\Gamma = \partial \Omega_1 \cap \partial \Omega_2$.  We
consider the following elliptic model problem: find $u: \Omega_0
\rightarrow \R$ such that
\begin{alignat}{3}
  -\Delta u_i  &= f_i &\quad &\text{in }  \Omega_i, \quad i = 1,2,
    \label{eq:strong_poisson_equation}
  \\
    \jump{\nabla u \cdot \vect{n}} &= 0 &\quad &\text{on } \Gamma,
    \label{eq:strong_cont_condition}
    \\
      \jump{u} &= 0 &\quad &\text{on } \Gamma,
      \label{eq:strong_flux_condition}
      \\
      u &= 0 &\quad &\text{on } \partial \OmegaD{},
      \label{eq:dirichlet_bc}
      \\
      \nabla u \cdot \vect{n} &= g &\quad &\text{on } \partial \OmegaN{}.
      \label{eq:neumann_bc}
\end{alignat}
We here use the notation $v_i = v|_{\Omega_i}$ for the restriction of
a function $v$ to the subdomain $\Omega_i$ for $i = 1, 2$, $\vect{n}$
is the unit normal to $\Gamma$ directed from $\Omega_2$ into
$\Omega_1$, and $\jump{v} = v_2 - v_1$ denotes the jump in a function
over the interface $\Gamma$. In addition, the boundary $\partial
\Omega_0$ is divided into two subdomains $\partial \OmegaD{}$ and $\partial
\OmegaN{}$ where Dirichlet and Neumann boundary conditions are
applied, respectively.

\subsection{Finite element formulation}
\label{ssec:finite_element_formulation}

\begin{figure}
  \includegraphics[height=0.35\textwidth]{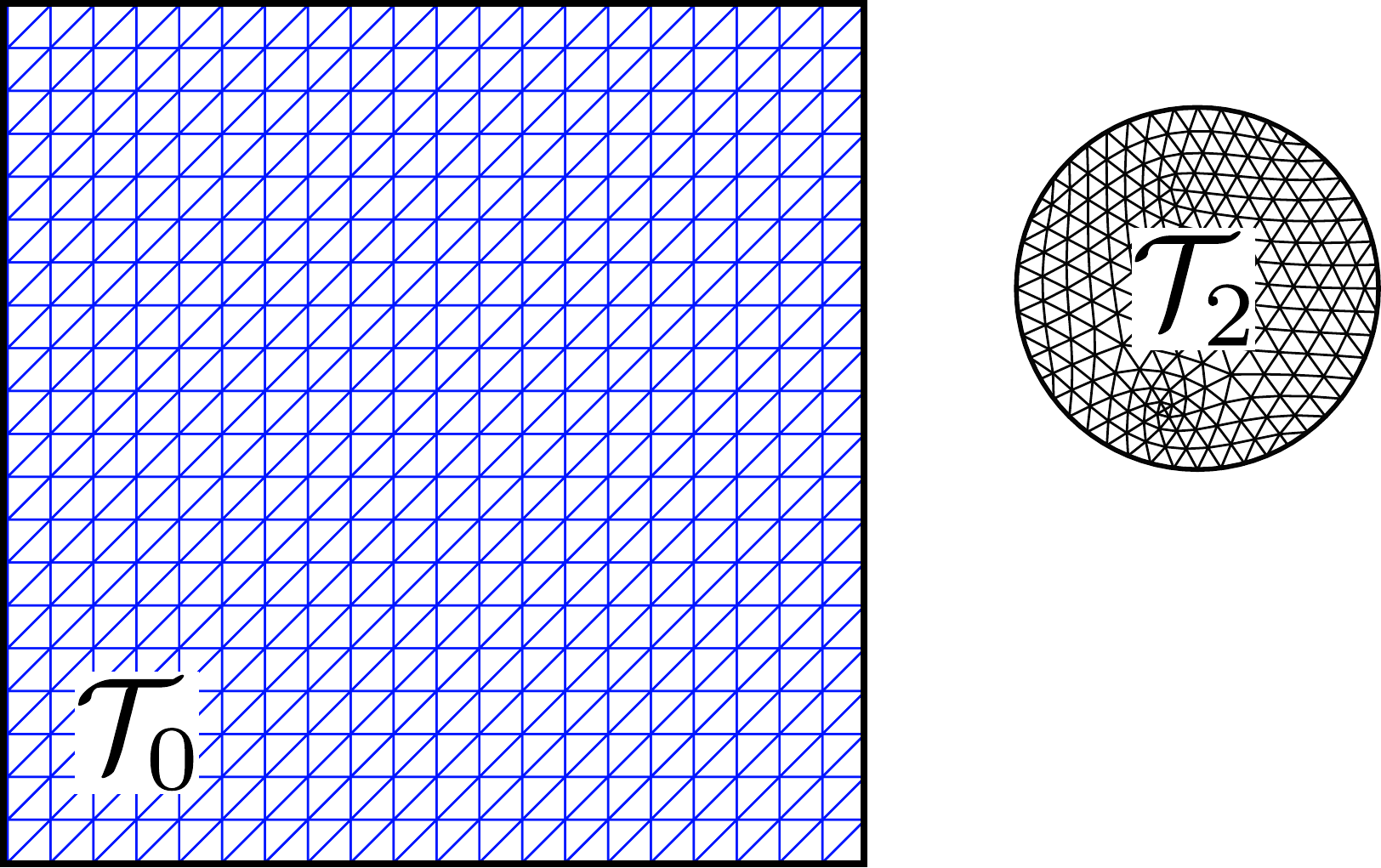}
  \hspace{0.05\textwidth}
  \includegraphics[height=0.35\textwidth]{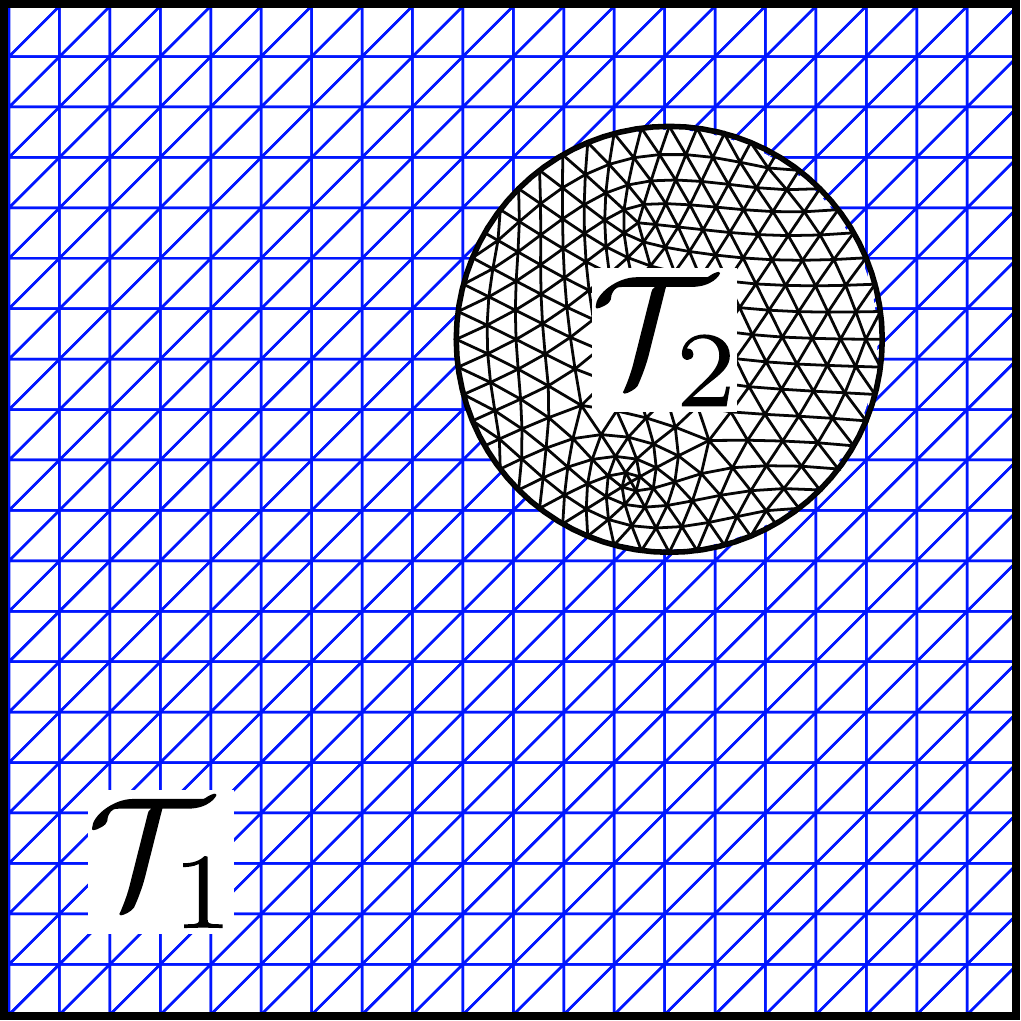}
  \caption{The mesh $\mesh_1$ for $\Omega_1 = (\Omega_0 \setminus
    \Omega_2)^{\circ}$ is constructed as the set of cells of $\mesh_0$
    not completely overlapped by the cells of $\mesh_2$.}
  \label{fig:overlapping_meshes}
\end{figure}

We consider a situation where a background mesh $\mesh_0$ is given for
$\Omega_0 = (\overline{\Omega_1 \cup \Omega_2})^{\circ}$ and another
mesh $\mesh_2$ is given for the overlapping domain $\Omega_2$ (see
Figure~\ref{fig:model_domain}). Both meshes are assumed to consist of
shape-regular tetrahedra $T$. The mesh $\mesh_1$ covering
$\Omega_1 = (\Omega_0 \setminus \Omega_2)^{\circ}$ is constructed by
\begin{equation}
  \mesh_1 = \{ T \cap \Omega_1 : T \in \mesh_0 \wedge |T \cap \Omega_1| > 0 \}.
  \label{eq:mesh1_def}
\end{equation}
Note that the $\mesh_1$ consists of both regular and cut elements
since the mesh $\mesh_0$ is not a conform tetrahedralization of the
subdomain $\Omega_1$ in general. For a cut element $T \cap \Omega$,
the degree of freedoms are defined to be the same as for the original
uncut element $T$.  This is illustrated in
Figure~\ref{fig:overlapping_meshes}.

We let ${V}_{h,i}$ denote the space of continuous piecewise fixed-order
polynomials on $\mesh_i$ that vanish on $\partial \Omega_i \cap
\partial \OmegaD{}$ for $i = 0, 2$. The space $V_{h,1}$ is constructed
as the restriction to $\Omega_1$ of functions in $V_{h,0}$; that is,
$V_{h,1} = \{ v|_{\Omega_1}: v \in V_{h,0}\}$.  We may then define a
finite element space for the whole domain $\Omega_0$ by
\begin{equation}
  {W}_h = {V}_{h,1} \bigoplus {V}_{h,2}.
  \label{eq:fem-space-definition}
\end{equation}
Nitsche' method proposed by~\citet{HansboHansboLarson2003} then takes
the form: find ${u}_h \in {W}_h$ such that
\begin{equation}
  \label{eq:varproblem}
  a({u}_h,{v}) = l({v}) \quad \foralls {v} \in {W}_h,
\end{equation}
where
\begin{align}
  \label{eq:poisson_nitsche_form}
  a(u,v)
  &=
  \sum_{i=1,2} \int_{\Omega_i} \nabla u \cdot \nabla{v} \dx
  \\ \nonumber
 &\qquad  - \underbrace{\int_{\Gamma}  \avg{\nabla u \cdot {\vect{n}}}
  \jump{{v}} \dS}_{\text{Consistency}}
  - \underbrace{\int_{\Gamma} \avg{\nabla v \cdot {\vect{n}}} \jump{u} \dS}_{\text{Symmetrization}}
  + \underbrace{\gamma \int_{\Gamma} h^{-1} \jump{{u}} \cdot \jump{{v}} \dS}_{\text{Penalty/Stabilization}},
  \\
  l(v) &= \int_\Omega f v \dx  +
          \int_{\partial \OmegaN{}} g v \ds.
  \label{eq:poisson_nitsche_form_rhs}
\end{align}
Here, $\gamma$ is a positive penalty
parameter and the average
$\avg{ \nabla v \cdot {\vect{n}}}$ is chosen to be the one-side
derivative
$\avg{ \nabla v \cdot {\vect{n}}} = \nabla v_2 \cdot \vect{n}$
according to \citet{HansboHansboLarson2003}. But as pointed out by the
authors, any convex combination of the normal derivatives leads to a
consistent formulation.

\subsection{Summary of theoretical results}

When the penalty parameter $\gamma$ in \eqref{eq:poisson_nitsche_form}
is chosen large enough, the form $a(\cdot, \cdot)$ is coercive on the
discrete space ${W}_h$ and one can derive optimal order \apriori{}
error estimates in both the energy norm and the $L^2$ norm for
polynomials of arbitrary degree $p$. The estimates take the form
\begin{gather}
  \left( \sum_{i=1}^2\| \nabla (u - u_h) \|^2_{\Omega_i} + h^{-1} \|
       [u_h]\|_\Gamma^2\right)^{\frac{1}{2}} \leq C h^p \| u
       \|_{\Omega_0,p+1},
       \label{eq:h1norm_error}
       \\
       \| u - u_h \| \leq C  h^{p+1} \| u \|_{\Omega_0,p+1}
       \label{eq:l2norm_error}.
\end{gather}
See~\citet{HansboHansboLarson2003} for details. Here, $\| \cdot
\|_{\Omega_0, s}$ denotes the standard Sobolev norm of order $s > 0$
on $\Omega_0$. \emph{A~posteriori} error estimates and adaptive
algorithms are also presented in~\citet{HansboHansboLarson2003}.

As observed by \citet{AreiasBelytschko2006}, Nitsche's method for an
arbitrary cutting interface as in
\citet{HansboHansbo2002,HansboHansbo2004} can be reinterpreted as a
particular instance of an extended finite element method. In contrast,
the presented formulation on overlapping meshes lacks such a
reinterpretation since \emph{two} unrelated meshes are involved and
therefore the discontinuity across the interface cannot be modeled by
the enrichment of degrees of freedom within a \emph{single} cell.
Nevertheless, both methods share some features, especially with regard
to their implementation.

\section{Techniques and algorithms}
\label{sec:implementation-challenges}

The main challenges in the implementation of Nitsche's method on
overlapping meshes arise from the geometric computations which are
necessary to assemble the discrete system associated with
\eqref{eq:varproblem}--\eqref{eq:poisson_nitsche_form_rhs}. Naturally,
these geometric computations are more involved in 3D than in 2D.
Similar challenges are encountered in related methods such as the
extended finite element method \citep{FriesBelytschko2010}, but also
for the simulation of contact mechanics \citep{YangLaursen2008}.
Different solutions have been proposed, see for example
\citet{SukumarMoesMoranEtAl2000} in the case of extended finite
elements methods, and \citet{YangLaursen2008} for contact problems.
Here, we take another approach which efficiently solves the issues in
the case of the Nitsche overlapping mesh method. In the following, we
will first discuss the challenges and their remedies in general terms,
and then return to the specific details of our implementation in
Section~\ref{sec:data-structures}.

\subsection{Assembly on overlapping meshes}
\label{ssec:main-challenges}

The main implementation challenges arise from the fact that the
interface $\Gamma$ can cross the overlapped mesh $\mesh_{0}$ in an
arbitrary manner, which has two consequences. First, the definition of
the finite element space \eqref{eq:fem-space-definition} involves the
restriction of the function space $V_{h,0}$, defined on the full
background mesh $\mesh_0$ of $\Omega_0$, to the domain
$\Omega_1$. The restriction results in non-standard element geometries
along the interface $\Gamma$ as it can be observed
from the definition~\ref{eq:mesh1_def} of $\mesh_1$.
Secondly, the weak imposition of the
interface conditions \eqref{eq:strong_cont_condition} and
\eqref{eq:strong_flux_condition} by the interface integrals in
\eqref{eq:poisson_nitsche_form} involves finite element spaces defined
on two unrelated meshes. Both make the assembly challenging to
implement, compared to a standard finite element method.

To better understand what kinds of challenges the Nitsche overlapping
mesh method adds, we briefly recall the general theme of the assembly
routine as it is realized nowadays in many finite element frameworks;
see, e.g., \citet{Logg2007,BastianHeimannMarnach2010}.  Nitsche's
method is closely related to the classical discontinuous Galerkin (DG)
methods, which makes it natural to depart from the assembly in the DG
case. A detailed description of finite element assembly for
discontinuous Galerkin methods can be found in
\citet{OelgaardLoggEtAl2008}. We assume a variational problem of the
following form: find $u_h \in V_h$ such that
\begin{equation}
  a(u_h,v) = l(v) \quad \foralls v \in \hat{V}_h,
  \label{eq:variational_problem}
\end{equation}
where $a$ and $l$ are bilinear and linear forms, respectively, and
$V_h$ and $\hat{V}_h$ are the discrete trial and test spaces,
respectively. For simplicity, we here make the assumption $V_h =
\hat{V}_h$.  The solution of the variational problem
\eqref{eq:variational_problem} may then be computed by solving the
linear system
\begin{equation}
  A U = b,
  \label{eq:linear_system}
\end{equation}
for $U \in \R^N$, with stiffness matrix $A_{IJ} = a(\phi_J,{\phi}_I)$
and load vector $b_I = l({\phi}_I)$. Here, $\{\phi_I\}_{I=1}^N$
denotes a basis for $V_h$. During assembly of the tensors $A$ and $b$,
one usually iterates over all cells $T \in \mesh$ and computes the
contributions from each cell as the cell tensors $A^T_{ij} =
a(\phi_j^T,\phi_i^T)$ and $b^T_i = l({\phi}_i)$ where
$\{\phi^T_i\}_{i=1}^n$ is the local finite element basis on $T$ (the
\emph{shape functions}). The element tensors are then scattered into
the global tensors $A$ and $b$ by adding the entries according to a
given local-to-global mapping $\iota_T: i \mapsto I$. One may
similarly add contributions from each facet~$F$ of the set of boundary
facets $\partial_e T$ (the set of exterior facets) and, for the
implementation of a DG method, from each facet $F$ of the set of
interior facets $\partial_i T$.  Algorithm~\ref{alg:standard_assembly}
summarizes the parts of the standard assembly algorithm relevant for
our discussion.

\begin{algorithm}
  \begin{tabbing}
  $A = 0$ \\
  \textbf{for each} $T \in \mesh$:\\
  \tab $\mathcal{I}(T) = \{1,\dots,\text{dim}(\hat{V}_h(T))\} \times
  \{1,\dots,\text{dim}(V_h(T))\}$\\
  \tab \textbf{for each} $(i,j) \in \mathcal{I}(T):$\\
  \tab \tab $A^T_{ij} = a(\phi^T_j,\hat{\phi}^T_i)$ \\
  \tab \textbf{for each} $(i,j) \in \mathcal{I}(T):$\\
  \tab \tab $A_{\hat{\iota}_T(i),\iota_T(j)} += A^T_{ij}$ \\
  \textbf{for each} $F \in \partial_e \mesh$: \\
  \tab $\mathcal{I}_e(F) = \{1,\dots,\text{dim}(\hat{V}_h^e(F))\} \times
  \{1,\dots,\text{dim}(V_h^e(F))\}$\\
  \tab \textbf{for each} $(i,j) \in \mathcal{I}_e(F)$:\\
  \tab \tab $A^T_{ij} = a(\phi^{F,e}_j,\hat{\phi}^{F,e}_i)$ \\
  \tab \textbf{for each} $(i,j) \in \mathcal{I}_e(F)$:\\
  \tab \tab $A_{\hat{\iota}_F^e(i),\iota^e_F(j)} += A^T_{ij}$ \\
  \textbf{for each} $F \in \partial_i \mesh$:\\
  \tab $\mathcal{I}_i(F) = \{1,\dots,\text{dim}(\hat{V}_h^i(F))\} \times
  \{1,\dots,\text{dim}(V_h^i(F))\}$\\
  \tab \textbf{for each} $(i,j) \in \mathcal{I}_i(F)$:\\
  \tab \tab $A^T_{ij} = a(\phi^{F,i}_j,\hat{\phi}^{F,i}_i)$ \\
  \tab \textbf{for each} $(i,j) \in \mathcal{I}_i(F)$:\\
  \tab \tab $A_{\hat{\iota}^i_F(i),\iota^i_F(j)} += A^T_{ij}$
  \end{tabbing}
  \label{alg:standard_assembly}
  \caption{Standard finite element assembly of a bilinear form
    $a(\cdot, \cdot)$. The local finite element spaces are denoted by
    $V_h(T), V_h(F)$ etc.}
\end{algorithm}

We now consider how the standard assembly algorithm must be modified
in the case of Nitsche's method on overlapping meshes. We first note
that the tessellation $\mesh_0$ of the background domain $\Omega_0$ may
be decomposed into three disjoint subsets:
\begin{equation}
  \label{eq:mesh-splitting}
  \mesh_0 = \mesh_{0,1} \cup \mesh_{0,2} \cup \mesh_{0,\Gamma},
\end{equation}
where $\mesh_{0,1} = \{ T \in \mesh_0 : T \subset \overline{\Omega}_1
\} $, $\mesh_{0,2} = \{ T \in \mesh_0 : T \subset \overline{\Omega}_2
\} $, and $\mesh_{0,\Gamma} = \{ T \in \mesh_0 : | T \cap \Omega_i |
> 0, \ i=1,2 \} $ denote the sets of \emph{not}, \emph{completely} and
\emph{partially} overlapped cells relative to $\Omega_2$,
respectively.

\begin{figure}
  \begin{center}
    \def\svgwidth{0.95\textwidth}
    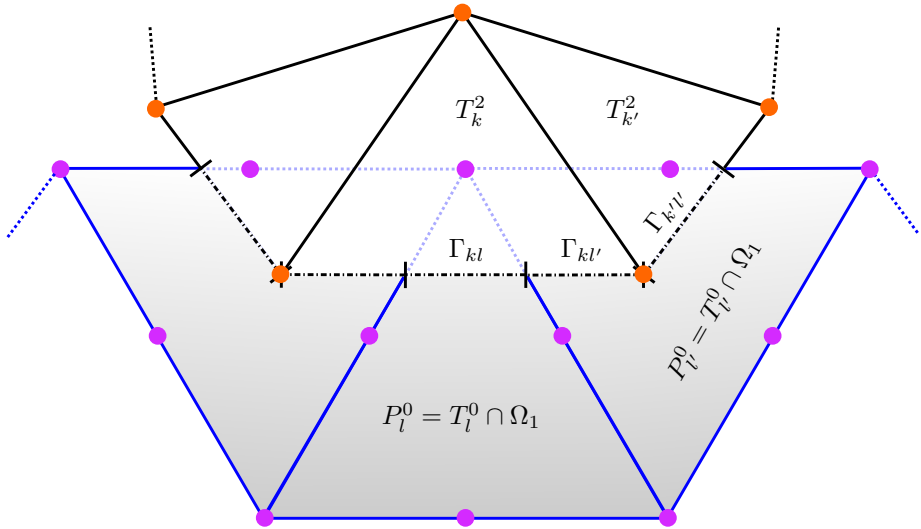
    \caption{The overlapping meshes in the intersection zone. The
      filled circles represent the degrees of freedom on each element
      (here for a piecewise linear approximation on $\mesh_2$ and
      piecewise quadratic on $\mesh_0$). The interface $\Gamma$ is
      partitioned so that each part intersects exactly one cell
      $T^2_k$ of the overlapping mesh $\mesh_2$ and one cell $T^0_l$
      of the background mesh $\mesh_0$.}
    \label{fig:intersection_areas}
  \end{center}
\end{figure}

Integrals over the cells of $\mesh_{0,1}$ can be assembled using a
standard assembly algorithm. Furthermore, integrals over the cells of
$\mesh_{0,2}$ need not be assembled (the corresponding contributions
will be assembled over $\mesh_2$). However, assembly must be carried
out over $\mesh_{0,\Gamma}$, the partially overlapped cells of
$\mesh_0$. For the model problem
\eqref{eq:varproblem}--\eqref{eq:poisson_nitsche_form_rhs}, this
requires the evaluation of integrals of the type
\begin{equation}
  \int_{P}  \nabla \phi^T_j \cdot \nabla \phi^T_i \dx
  \, \text{ and } \,
  \int_{P} f \phi^T_i \dx,
  \label{eq:cut_cell_tensor}
\end{equation}
on \emph{cut} elements (polyhedra) $P = T \cap \Omega_1$ where $T \in
\mesh_0$. Examining
\eqref{eq:varproblem}--\eqref{eq:poisson_nitsche_form_rhs},
we further note that we must assemble the terms
\begin{equation}
  - \int_{\Gamma} \avg{\nabla \phi_j \cdot \vect{n}} \jump{\phi_i} \dS
  - \int_{\Gamma} \avg{\nabla \phi_i \cdot \vect{n}} \jump{\phi_j} \dS
  + \gamma \int_{\Gamma} h^{-1} \jump{\phi_j} \cdot \jump{\phi_i} \dS.
  \label{eq:interface_tensor}
\end{equation}
This poses an additional challenge, since the integrands involve
products of trial and test functions defined on \emph{different}
meshes. Furthermore, the interface $\Gamma$ consists of a subset of
the boundary facets $\partial_e \mesh_2$ of $\mesh_2$, but each such
facet may intersect several cells of $\mesh_0$. We therefore partition
each facet on $\Gamma$ into a set of polygons $\{\Gamma_{kl}\}$ such
that each polygon $\Gamma_{kl}$ intersects exactly one cell $T^2_k$ of
the overlapping mesh $\mesh_2$ and one cell $T^0_l$ of the background
mesh $\mesh_0$. This is illustrated in
Figure~\ref{fig:intersection_areas}. Assembly may then be carried out
by summing the contributions from each polygon $\Gamma_{kl}$.

In summary, we identify the following main challenges in the
implementation of Nitsche's method on overlapping meshes:
\begin{enumerate}
\item
  \emph{collision detection}: to determine which cells are involved in
  the intersection between the two meshes;
\item
  \emph{mesh intersection}: to compute the cut cells of the background
  mesh $\mesh_0$ represented by the polyhedra $\{P^0_l\}$ and the
  intersection interface represented by the polygons
  $\{\Gamma_{kl}\}$;
\item
  \emph{integration on complex polyhedra}: to compute integrals on the
  polyhedra $\{P^0_l\}$ (and the polygons $\{\Gamma_{kl}\}$).
\end{enumerate}
In the following, we discuss these challenges in some detail and
introduce concepts, data structures and algorithms to handle them in
an efficient manner. We emphasize that the proposed solutions are not
limited to the implementation of Nitsche-type methods, but may also be
used for the implementation of related overlapping mesh methods.

\subsection{Collision detection}
\label{ssec:collision-detection}

Topological relations between entities of a single mesh can be
described by concepts of connectivity or mesh incidence as presented
in \citet{Logg2009a} and \citet{BastianDroskeEngwer2004}. To represent
the topological relation between two overlapping (colliding) meshes,
we enrich the notation of \citet{Logg2009a} by the concepts of
\emph{collision relations} and \emph{collision maps}. The collision
relation $\mesh_0 \leftrightarrow \partial \mesh_2$ between $\mesh_0$
and $\partial \mesh_2$ is defined by
\begin{equation}
  \mesh_0 \leftrightarrow \partial \mesh_2 = \{ (i,j) : T_i \cap F_j
  \neq \emptyset \wedge T_i \in \mesh_0 \wedge F_j \in \partial \mesh_2
  \}.
  \label{eq:collision-relation}
\end{equation}
The collision relation lists pairs of indices of all intersecting
cells of the background mesh $\mesh_0$ with boundary facets of the
overlapping mesh $\mesh_2$. Furthermore, each index of an intersected
entity is mapped to the set of indices of intersecting entities via a
pair of collision maps:
\begin{align}
  \label{eq:collision-map-1}
  (\mesh_0 \to \partial \mesh_2)(i) =  \{ j : (i,j) \in \mesh_0
  \leftrightarrow \partial \mesh_2 \}, \\
  \label{eq:collision-map-2}
  (\partial \mesh_2 \to \mesh_0)(j) =  \{ i : (i,j) \in \mesh_0
  \leftrightarrow \partial \mesh_2 \}.
\end{align}
We note that the collision maps $\mesh_0 \to \partial \mesh_2$ and
$\partial \mesh_2 \to \mesh_0$ can be computed from the
collision relation $\mesh_0 \leftrightarrow \partial \mesh_2$.

A naive approach to computing the collision relation between two
meshes $\mesh_0$ and $\mesh_2$ would be to intersect each cell of
$\mesh_0$ with each cell of $\mesh_2$, resulting in an
$\mathcal{O}(|\mesh_0| \cdot |\mesh_2|)$-complexity, which is not
feasible for large meshes. However, efficient algorithms and data
structures which reduce the complexity dramatically have been
developed in the fields of computer science, computational geometry,
and computer graphics. The task of determining whether two objects
collide arises naturally in the rendering of a computer game and the
objects in question are often represented by meshes. In this paper, we
limit the description to those techniques used in our work and refer
to the books \citet{ericson2005real,Schneider2003}, and
\citet{Akenine-MollerHainesHoffman2008} for a broader overview. To
efficiently find all pairs of the collision relation $\mesh_0
\leftrightarrow \partial \mesh_2$, two important concepts are used:
(i) fast intersection tests for pairs of simple geometric objects and
(ii) spatial data structures for large sets of objects to accelerate
collision queries.

\subsubsection{Fast intersection detection}

\begin{figure}
  \includegraphics[width=0.40\textwidth]{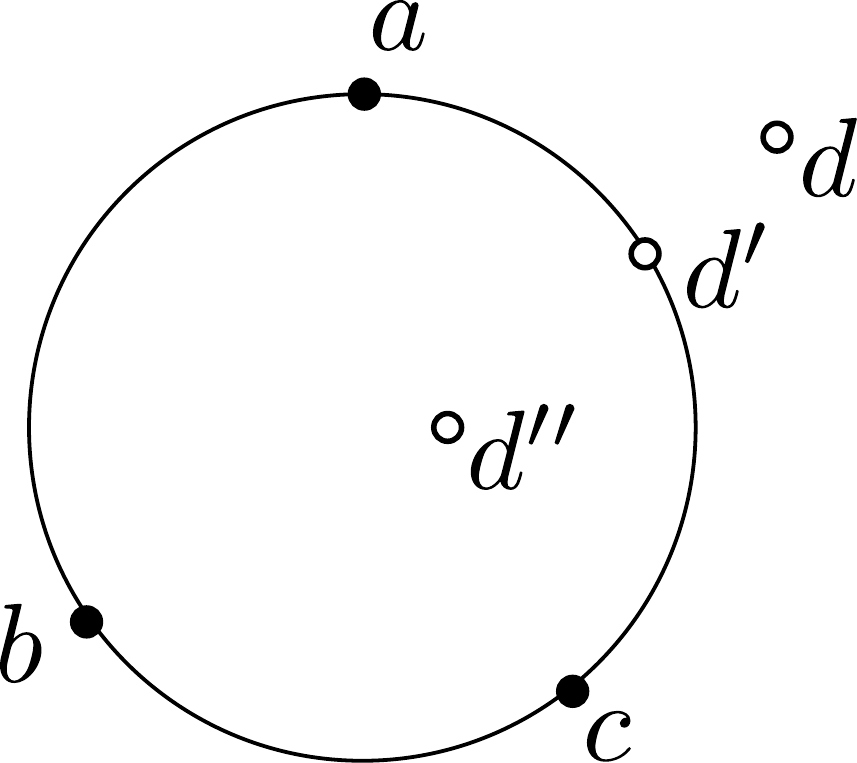}
  \hspace{1.0cm}
  \includegraphics[width=0.50\textwidth]{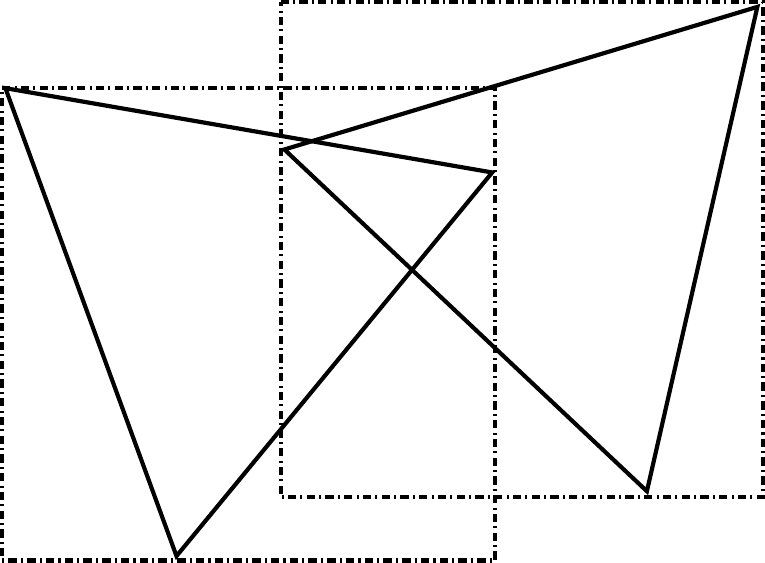}
  \caption{(Left) A typical example of a geometric predicate is the
    incircle test, which tells whether a given point $d$ lies on the
    inside, on the outside, or on a circle defined by three points
    $a,\ b,\ c$. (Right) Axis aligned bounding boxes (AABB) of two
    triangles.}
  \label{fig:incirlce-test-and-bbox}
\end{figure}

Fast intersection detection for simple geometric entities such as
triangles and tetrahedra is based on so-called \emph{geometric
  predicates}. Geometric predicates are tests which determine whether
two geometric entities do intersect (collide) without computing the
actual intersection. See Figure~\ref{fig:incirlce-test-and-bbox} for
an example. This saves work since the actual intersection does not
need to be computed when the geometric predicate is false.

\subsubsection{Spatial data structures}

Spatial data structures provide a hierarchical ordering of geometric
objects, which allows quick traversal of large sets of geometric
objects as part of a collision test. There exist two principal
approaches based on either space subdivisions or bounding volume
hierarchies, which can both be realized as tree-like data
structures. As the name indicates, subdivision approaches rely on some
sort of geometric subdivision of the \emph{entire} space embedding the
structures of interest, in our case a finite element mesh. Typical
examples are binary space partitions (BSP) trees, quadtrees, and
octrees. Since the embedding space is subdivided, the leaf of a tree
will usually not represent a single mesh entity and therefore further
selection procedures are required to determine the actual
intersections.

\begin{figure}
  \includegraphics[width=0.48\textwidth]{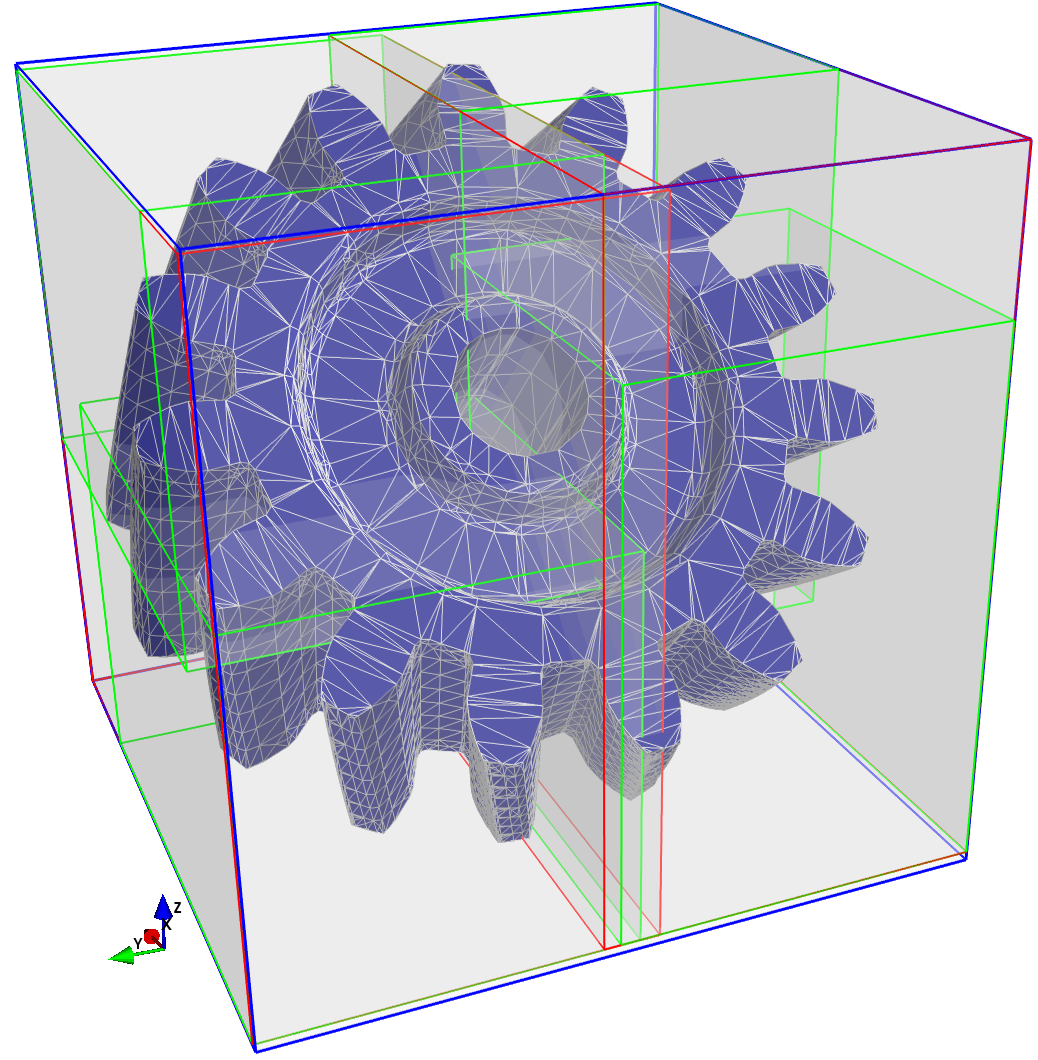}
  \includegraphics[width=0.48\textwidth]{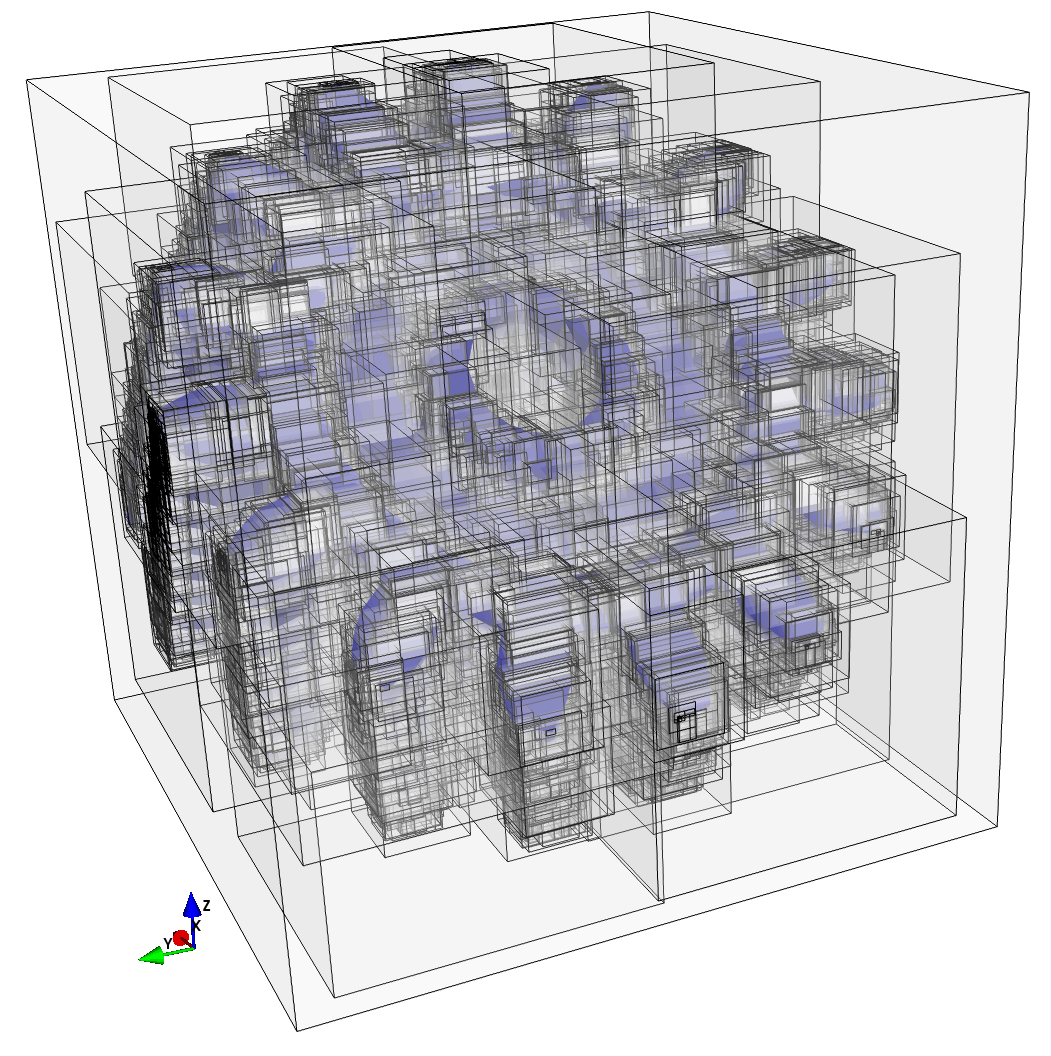}
  \caption{Bounding volume tree of gear mesh. (Left) Root bounding box
    (blue) and the next two generations (red and green).  (Right)
    Every fourth level of the complete bounding box tree.}
  \label{fig:bbox_tree}
\end{figure}

In contrast, a bounding volume hierarchy is a tree which is built from
bounding volumes; that is, simple geometric shapes containing the objects
to be tested. Examples of such data structures are axis aligned
bounding boxes (AABB, see Figures~\ref{fig:incirlce-test-and-bbox}
and~\ref{fig:bbox_tree}), oriented bounding boxes (OBB) and so-called
k-DOPs, discrete orientation polytopes described by
$k$~hyperplanes. The purpose of using bounding volumes is twofold: (i)
testing for collision of bounding volumes is cheaper than testing for
collision of the bounded objects, and (ii) the simple geometry of the
bounding volumes means that they can be stored efficiently in a
hierarchical manner.

\subsubsection{Building the collision map}

Bounding volume trees accelerate asymmetric collision queries between
a tree embedding a large set of objects and a single simple object
like a tetrahedron. The concept becomes even more powerful when
intersection tests between two large sets of geometric primitives are
desired, as in our case of two meshes. Then a hierarchical traversal
of both trees greatly improves the $\mathcal{O}(|\mesh_0| \cdot
|\mesh_2|)$-complexity of the naive approach.
Algorithm~\ref{alg:tree_traversal} describes how to employ a
hierarchical traversal to compute the collision map.

\begin{algorithm}
  \begin{tabbing}
    compute\_collisions($A$, $B$):\\
    \tab \textbf{if } $A \cap B = \emptyset$:\\
    \tab \tab \textbf{return} \\
    \tab \textbf{else if} is\_leaf($A$) $\wedge$ is\_leaf($B$): \\
    \tab \tab \textbf{if } $T_A \cap T_B \neq \emptyset$:\\
    \tab \tab \tab  $(i, j) =
    (\mathrm{index}(T_A), \mathrm{index}(F_B))$ \\
    \tab \tab \tab $\mesh_0 \leftrightarrow \partial \mesh_2 :=
    \mesh_0 \leftrightarrow \partial \mesh_2 \cup (i, j)$ \\
    \tab \textbf{else if} descend\_a($A$, $B$): \\
    \tab \tab \textbf{for} $a$ $\in$ children($A$): \\
    \tab \tab \tab compute\_collisions($a$, $B$) \\
    \tab \textbf{else}: \\
    \tab \tab \textbf{for} $b$ $\in$ children($B$): \\
    \tab \tab \tab compute\_collisions($A$, $b$) \\
    \\
    descend\_a($A$, $B$): \\
    \tab \textbf{return} is\_leaf($B$) $\vee \,
    (\neg$is\_leaf($A$) $\wedge \, |A| \ge |B|)$
  \end{tabbing}
  \caption{Traversal of two bounding volume trees. Starting from the
    root of both trees, pairs of nodes are recursively tested for
    intersection. Only if the bounding volumes $A$ and $B$ of two
    nodes overlap may the children possibly overlap. A so-called
    descend rule is applied to decide with which of the two nodes to
    proceed. A popular and effective rule is to choose the one with
    the larger volume \citep{ericson2005real} since it gives the
    largest volume reduction for subsequent bounding volumes. If two
    leaf bounding volumes are reached, then the mesh entities bounded
    by them are tested for intersection.}
  \label{alg:tree_traversal}
\end{algorithm}

After the completion of the traversal according to
Algorithm~\ref{alg:tree_traversal}, which identifies all cells of
$\mesh_0$ that intersect the \emph{boundary} of $\mesh_2$, it remains
to check whether the remaining cells (which do not intersect the
boundary of $\mesh_2$) are either completely overlapped or not
overlapped by $\mesh_2$. To check this, it is enough to take a single
point contained in each cell and check whether it is in $\Omega_2$. To
avoid building a third AABB tree for the entire mesh $\mesh_2$ and to
take advantage of the already built tree for $\partial \mesh_2$, one
can use a method known as \emph{ray-shooting}; see
\citet{Akenine-MollerHainesHoffman2008}. The idea is simple: counting
how often a ray starting at the point in question intersects the
surface $\partial \mesh_2$ tells whether it is inside (odd number of
intersections) or outside (even number of intersections). To
efficiently find all ray--surface intersections, one may reuse the
same AABB tree for $\partial \mesh_2$ that was used in
Algorithm~\ref{alg:tree_traversal}.

To summarize, the collision map and classification of cells according
to the splitting~\eqref{eq:mesh-splitting} can be found efficiently by
utilizing AABB trees, fast traversal of pairs of AABB hierarchies, and
ray shooting techniques.

\subsection{Mesh intersection}

\begin{figure}
  \includegraphics[width=0.50\textwidth]{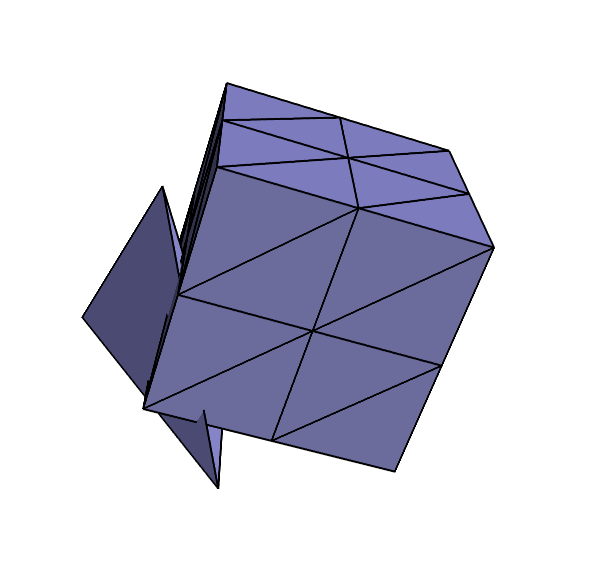}
  \includegraphics[width=0.45\textwidth]{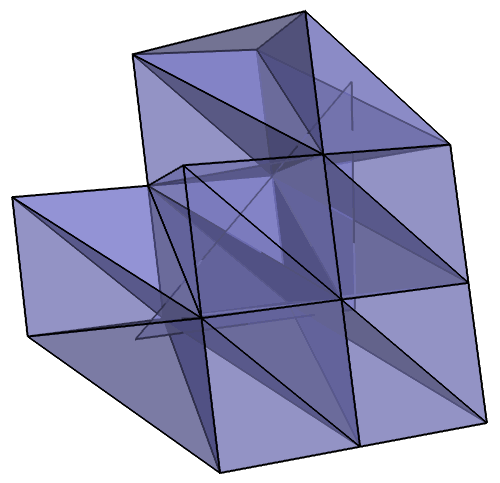}
  \\
  \includegraphics[width=0.45\textwidth]{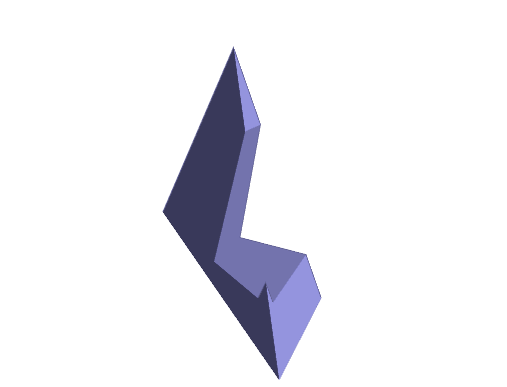}
  \hspace{1.2cm}
  \includegraphics[width=0.24\textwidth]{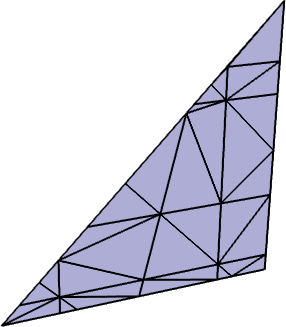}
  \caption{(Left top) A tetrahedron cell $T^0_l$ from the background
    mesh $\mesh_0$ overlapped by a tetrahedral mesh $\mesh_2$ of a
    cube. (Bottom left) The resulting cut cell (polyhedron)
    $P^0_l$. (Top right) An interface facet $F \in \partial \mesh_2$
    of the overlapping mesh $\mesh_2$ intersected by a number of
    tetrahedra from the background mesh $\mesh_0$.  (Bottom left) The
    resulting facet decomposition (triangulation) $\{\Gamma_{kl}\}$.}
  \label{fig:intersected-cells-and-facets}
\end{figure}

The next step is to compute the cut cells (polyhedra) $\{P^0_l\}$ and
the interface decomposition $\{\Gamma_{kl}\}$; see
Figure~\ref{fig:intersection_areas}. This computation may be phrased
in terms of so-called \emph{boolean operations} which are widely used
in CAD systems to build complex geometries by performing boolean
operations between primitives from a finite set of geometries.
Algorithm~\ref{alg:compute_geometry} summarizes the use of boolean
operations to compute $\{P^0_l\}$ and $\{\Gamma_{kl}\}$. The resulting
geometric objects are depicted in
Figure~\ref{fig:intersected-cells-and-facets}.

\begin{algorithm}
  \begin{tabbing}
  \textbf{for each} $T \in \mesh_0$ intersected by $\Gamma$:\\
  \tab $l = \mathrm{index}(T)$ \\
  \tab $P^0_l := T$ \\
  \tab $S := \emptyset$ \\
  \tab \textbf{for each} $F \in (\mesh_0 \to \partial\mesh_2)(l)$: \\
  \tab \tab $S := S \cup F$ \\
  \tab $P^0_l :=  T^+_S $ \\
  \textbf{for each} $F \in \partial \mesh_2$:\\
  \tab $j = \mathrm{index}(F)$\\
  \tab $k = \mathrm{index}(T_F)$\\
  \tab \textbf{for each} $T \in (\partial \mesh_2 \to \mesh_0)(j)$:\\
  \tab \tab $l = \mathrm{index}(T)$ \\
  \tab \tab $\Gamma_{kl} := F \cap T$
  \end{tabbing}
  \caption{Mesh intersection. Note that the surface $S$ defined by the
  union of all facets $F$ intersecting the tetrahedron $T$ divides $T$
  into two parts $T^+_S$ and $T^-_S$ according to orientation of $S$.}
  \label{alg:compute_geometry}
\end{algorithm}

The boolean operations are completely delegated to the computational
geometry library \emp{GTS} \citep{gts} which uses a so-called
\emph{ear-clipping} algorithm to compute the surface tesselation of
the intersection objects. The implemented ear clipping algorithms
(first introduced in \citet{Fournier1984}) are known to have
$\mathcal{O}(n \log(n) )$ complexity, where $n$ is the number of
vertices of the 3D-polygon. This number in turn ($n$) depends on the
mesh quality parameters such as the smallest and widest element angles
and will in practice be bounded by a constant.

\subsection{Integration on complex polyhedra}
\label{ssec:integration_complex_polyhedron}

The last required step in implementing the overlapping mesh method is
to compute the integrals \eqref{eq:cut_cell_tensor} and
\eqref{eq:interface_tensor} on the cut cells $\{P^0_l\}$ and on the
interface decomposition $\{\Gamma_{kl}\}$, respectively. A widely used
technique in the implementation of the extended finite element
method \citep{es1999finite,Daux2000,SukumarMoesMoranEtAl2000,Belytschko2001}
is to decompose a cell into subcells which are aligned with the
interface. The interface is either approximated by a plane within the
cell \citep{Fries2008,Fries2009} or completely recovered
\citep{Daux2000,SukumarMoesMoranEtAl2000}, even in the case of higher order
interfaces \citep{MayerGerstenbergerWall2009}. Inside each tetrahedron
of the subtetrahedralization, one may then use a standard integration
scheme. However, subtetrahedralization of an arbitrary polyhedron is
in general a quite challenging problem~\citep{RuppertSeidel1989}, and
its existence cannot even be guaranteed without adding additional
vertices \citep{ChazellePalios1990}, in contrast to the two-dimensional
case.

Methods for integration over arbitrary polygonal domains without the
use of a subtriangulation have recently been presented in
\citet{Natarajan2009} and incorporated into the extended finite element
method with discontinuous enrichments in \citet{Natarajan2010}. Using
Schwarz--Christoffel conformal mappings as the fundamental tool, the
technique is strongly bounded to two space dimensions and hard to
generalize to three space dimensions.

We here propose an alternative approach, which is based on a boundary
representation of the integrals. Despite the fact that this technique
has been known for a long time, see for instance \cite{Lee1982}, it
seems to be largely unknown within the finite element community. Our
implementation is based on the efficient realization described in
\citet{Mirtich1996}. This technique may be easily generalized to the
computation of the general moment integral $I_{\alpha}(P)$ over a
polyhedron $P$, defined by
\begin{equation}
  I_{\alpha}(P) = \int_{P} x^\alpha\dx = \int_P x_1^{\alpha_1}\cdots
  x_d^{\alpha_d} \dx,
  \label{eq:moment_equation}
\end{equation}
where $\vect{\alpha} = (\alpha_1, \ldots, \alpha_d)$ denotes a
multi-index of length $d$.

The integral $I_{\alpha}(P)$ is computed in three steps.  The first
step is to interpret the integrand as the divergence of a polynomial
vector field and to rewrite \eqref{eq:moment_equation} as surface
integral:
\begin{equation}
  \int_{P} x^{\alpha} \dx = \int_{P} \nabla \cdot
  \sum_{i=1}^d \dfrac{x^{\alpha+\vect{e}_i}}{d(\alpha_i + 1)} \vect{e}_i \dx
  = \sum_{F \in \partial P} \sum_{i=1}^d
  \vect{n}_{F} \cdot \vect{e}_i
  \int_{F} \dfrac{x^{\alpha+\vect{e}_i}}{d(\alpha_i + 1)} \ds.
  \label{eq:moment_integral_reduction}
\end{equation}
Here, $\vect{e}_i = (0, \ldots, 0, 1, 0, \ldots, 0)$ denotes the $i$th
unit vector.

Secondly, the plane equation $a x + b y + c z = d$ for each facet $F$
allows the construction of a projection map $Z = h(X,Y)$, where
$X,Y,Z$ is some positive, orientation-preserving permutation of
$x,y,z$. Using the parametrization $h$, one may rewrite the integrals
of type $\int_{F} x^{\beta} \ds$
in~\eqref{eq:moment_integral_reduction} as integrals in the
$XY$-plane:
\begin{equation}
  \int_{F} x^{\beta} \ds = \dfrac{1}{|{n_{Z}}|} \int_{h^{-1}(F)}
  x(X,Y,h(X,Y))^{\beta} \dX \dY.
  \label{eq:integral_with_projection}
\end{equation}
Here, $n_Z$ denotes the $Z$-component of the normal vector $\vect{n}_F$.

The third step consists of using a parametrization of $\partial
h^{-1}(F)$ and Green's theorem in the plane to rewrite
\eqref{eq:integral_with_projection} as a sum of line integrals.
Together, the three steps reduce the evaluation of the moment integral
\eqref{eq:moment_equation} to the evaluation of a set of
one-dimensional integrals with polynomial integrands. We note that the
algorithm for the calculation of each moment integral has only a
$\mathcal{O}(\#F)$ complexity.

The moment integrals can be used in several ways to finally compute
the integrals \eqref{eq:cut_cell_tensor} and
\eqref{eq:interface_tensor}. First, one may compute the volume and the
barycenter as the zeroth and first moments, respectively, which
immediately provides a quadrature rule for exact integration of
polynomials of degree one. Unfortunately, higher order moments may not
be directly reinterpreted as quadrature rules. Furthermore, as the
construction of quadrature rules for higher degree polynomials can be
quite involved \citep{Cools1997}, it can be advantageous to avoid
run-time generation of quadrature rules. We therefore propose to use
the moment integrals $I_{\alpha}(P)$ directly by first interpolating
the integrand onto a monomial basis and then summing the
contributions:
\begin{equation}
  \int_{P} f \dx
  = \int_{P} \sum_{\alpha} f_{\alpha} x^{\alpha} \dx
  = \sum_{\alpha} f_{\alpha} \int_{P} x^{\alpha} \dx
  = \sum_{\alpha} f_{\alpha} I_{\alpha}(P).
  \label{eq:integral-interpolation}
\end{equation}
Ill-conditioning may be avoided for higher-order expansions by
replacing the monomial basis by Bernstein polynomials. In the present
work, linear Lagrange elements have been used throughout and so simple
barycenter quadrature suffices.

\section{Implementation and data structures}
\label{sec:data-structures}

We now discuss some of the data structures and classes which reflect
the abstract concepts and algorithms described in
Section~\ref{sec:implementation-challenges}.  For some of these
algorithms, we rely on existing implementations as part of the
computational geometry libraries CGAL \citep{cgal} and GTS
\citep{gts}, while other algorithms have been realized in the finite
element library DOLFIN \citep{LoggWells2010a,LoggWellsEtAl2011} as
part of this work. Specialized algorithms for the Nitsche overlapping
mesh method have been implemented as part of the extension library
DOLFIN-OLM which is built on top of DOLFIN. The code is
free/open-source, licensed under the LGPLv3, and available at
\emp{http://launchpad.net/dolfin-olm}.

\subsection{The finite element library DOLFIN}
\label{sec:fenics-interface}

Our implementation of Nitsche's method on overlapping meshes is based
on the finite element library DOLFIN which is part of the FEniCS
project~\citep{LoggMardalEtAl2011,Logg2007} for automated scientific
computing. The main feature of FEniCS is the automated treatment of
finite element variational problems, based on automated generation of
highly efficient C++ code from abstract high-level descriptions of
finite element variational problems expressed in near-mathematical
notation~\citep{KirbyLogg2006,Alnaes2011a}. This is combined with
built-in tools for working with efficient representations of
computational meshes~\citep{Logg2009a} and wrappers for
high-performance linear algebra libraries like
PETSc~\citep{petsc-web-page,petsc-user-ref,petsc-efficient} and
Trilinos~\citep{trilinos}.

For this work, we have integrated the computational geometry libraries
CGAL~\citep{cgal} and GTS~\citep{gts} with DOLFIN. While CGAL and GTS
provide a large part of the functionality needed to compute
intersections between tetrahedra, the integration scheme was realized
through an adapted version of Mirtich's
code~\citep{Mirtich1996}. Furthermore, the assembly routine of DOLFIN
was extended to handle integration over cut cells and meshes.
Currently, the code for computing local interface integrals has to be
implemented manually by the user, but we plan to extend FEniCS, in
particular DOLFIN and the form compiler
FFC~\citep{KirbyLogg2006,LoggOelgaardEtAl2011a}, to provide a full
automation of Nitsche's method where a user only needs to supply the
abstract variational problem~\eqref{eq:varproblem}.

In the remaining subsections, we present the class abstractions of the
algorithms and data structures described in
Section~\ref{sec:implementation-challenges}. These new classes
introduced in DOLFIN-OLM allow the implementation of the Nitsche
assembly algorithm in a descriptive and concise manner, as illustrated
in Figure~\ref{code:interface-assembly}.

\begin{figure}
  \begin{center}
    \codesize
    \begin{code}
void NitscheAssembler::assemble_interface(GenericTensor& A,
					  const NitscheForm& a_nit,
					  UFC& ufc_0,
					  UFC& ufc_1)
{
  ...
  for (CutFacetPartIterator cut_facet(overlapping_meshes);
                            !cut_facet.end(); ++cut_facet)
  {
    // Update quadrature cache to current facet part and return quadrature rule
    const QuadratureRule& quadrature =
      a_nit.interface_domain_quadrature_cache()[cut_facet->index()];
    ...

    // Get cells incident with this part of the cut cell
    std::pair<const Cell, const Cell> cells = cut_facet->adjacent_cells();
    const Cell& cell0 = cells.first;
    const Cell& cell1 = cells.second;

    // Update ufc forms and update interface local dimensions
    uint local_facet1 = cell1.index(cut_facet->entire_facet());
    ufc_0.update(cell0);
    ufc_1.update(cell1);

    // Tabulate dofs for each dimension on interface facet part
    ...

    // Tabulate interface tensor.
    a_nit.tabulate_interface_tensor(a_nit.interface_A.get(),
				    ufc_0.cell,
				    ufc_1.cell,
				    local_facet1,
				    quadrature.size(),
				    quadrature.points(),
				    quadrature.weights());

    // Insert matrix
    A.add(a_nit.interface_A.get(), interface_dofs);
  }
}
    \end{code}
  \end{center}
  \caption{C++ implementation of the assembly of the so-called
    interface tensor accounting for the coupling between the two
    meshes $\mesh_1$ and $\mesh_2$. A \emp{CutFacetPartIterator}
    provides iteration over all interface facet parts, represented by
    \emp{CutFacetPart} and stemming from the facet decomposition
    $\{\Gamma_{kl}\}$. Each \emp{CutFacetPart} is associated with
    exactly one cell in each mesh, which can be accessed via the
    \emp{adjacent\_cells} member function. The overall design stresses
    the similarities to the assembly of facet contributions in the
    standard DG method.}
  \label{code:interface-assembly}
\end{figure}

\subsection{The class \emp{AABBTree}}

The search data structure \emp{AABBTree} has been added to the DOLFIN
library. The implementation is based on the computational geometry
library \emp{CGAL} \citep{cgal}. Basic search queries such as finding
one or all cells intersecting a given entity or distance computation
are exposed via an \emp{IntersectionOperator} class. DOLFIN-OLM
complements that functionality with providing a GTS-based AABB tree
\citep{gts} to allow traversal of two bounding box trees as described
in Algorithm~\ref{alg:tree_traversal}.

\subsection{The class \emp{OverlappingMeshes}}

The class \emp{OverlappingMeshes}, provided as part of DOLFIN-OLM, is
a key component in our realization of the overlapping mesh method. It
mainly computes additional topological and geometric information to
describe the overlap of the two meshes $\mesh_0$ and $\mesh_2$ and
provides access to this information, in particular the collision
relation $\mesh_0 \leftrightarrow \partial \mesh_2$ and the collision
maps $\mesh_0 \to \partial \mesh_2$ and $\partial \mesh_2 \to \partial
\mesh_0$ described in Section~\ref{ssec:collision-detection}.

A \emp{MeshFunction}, as introduced in \citet{Logg2009a}, describes the
splitting~\eqref{eq:mesh-splitting} of the mesh $\mesh_0$ by assigning
different integer values to cells in the not overlapped part
$\mesh_{0,1}$, the completely overlapped part $\mesh_{0,2}$, and the
partially overlapped part $\mesh_{0,\Gamma}$, respectively. In the
same way, the boundary facets of the overlapping mesh $\mesh_2$ can be
marked if the overlapping domain $\Omega_2$ is not completely
contained in background domain $\Omega_0$.

\subsection{Mesh iterators for overlapping meshes}

Mesh iterators have been advocated by \citet{Berti2002,Berti2006} and
\citet{Logg2009a} and used among others by
\citet{BastianHeimannMarnach2010} and \citet{Botsch2002} as an
important abstraction concept in mesh implementations and FEM
frameworks \citep{LoggWells2010a, BastianHeimannMarnach2010}. The
iterator concept allows to access and iterate over mesh entities such
as vertices, edges, facets and cells without knowing the details of
the underlying mesh implementation. We have followed the same ideas in
the case of overlapping meshes. The two pairs of classes
\emp{CutCell}, \emp{CutCellIterator} and \emp{CutFacetPart},
\emp{CutFacetPartIterator} provide an interface to the intersected
cells in $\mesh_{0,\Gamma}$ and the interface facet partition
$\{\Gamma_{kl}\}$, respectively. The \emp{FacetPart} class and the
corresponding iterator class mimic the original interface in DOLFIN by
giving access to the two adjacent cells in the overlapping and the
overlapped meshes, respectively. The interface assembly in
Figure~\ref{code:interface-assembly} presents a important use case.

Similar iterator concepts have been used in
\citet{BastianBuseSander2010} where a general infrastructure to couple
grids interfaced by the DUNE grid framework is presented.  The
\emp{CutFacetPartIterator} corresponds to the
\emp{Remote\-Intersection\-Iterator} described in
\citet{BastianBuseSander2010} specialized to our case of Nitsche's
method on overlapping meshes. Moreover, the \emp{CutCellIterator}
introduced here can be interpreted as an instance of
\emp{DomainIntersectionIterator} in \citet{BastianBuseSander2010}.

\subsection{The class \emp{Quadrature} and \emp{QuadratureRuleCache}}

The DOLFIN class \emp{Quadrature} is a lightweight base class which
only computes and stores quadrature data such as quadrature points,
weights and order for a given polyhedron at run-time. Since the
integration order depends on the underlying finite element scheme, the
actual computation of the points and weights is meant to be
implemented in subclass constructors. The class
\emp{BarycenterQuadrature} is such an instance of a subclass which
computes a quadrature rule of order~$2$ for a given polyhedron based
on the algorithm outlined in
Section~\ref{ssec:integration_complex_polyhedron}. For each
intersected cell or facet part, a \emp{QuadratureRule} object can be
stored in a~\emp{QuadratureRuleCache} instance to save geometry and
quadrature rule recomputation if several integrations have to be
performed on the same intersected entity, as is the case for the
computation of the stiffness matrix and load vector in the Nitsche
overlapping mesh method.

\subsection{Forms and assembly}

The DOLFIN \emp{Form} form class represents the mathematical concept
of a finite element variational form. This class has been extended, as
part of DOLFIN-OLM, to reflect the domain decomposition character of
the overlapping mesh method. A so-called \emp{NitscheForm} class holds
the description of the variational problem on each part $\Omega_1$ and
$\Omega_2$ of the domain. The coupling between the two forms is
accomplished through a member function
\emp{tabulate\_interface\_tensor}, which computes the local interface
tensor corresponding to \eqref{eq:interface_tensor}. This is in
addition to the standard \emp{tabulate\_tensor} functions defined in
the UFC code generation interface \citep{AlnaesLoggEtAl2009a} for cell
integrals, exterior facet integrals, and interior facet integrals
(cf. Algorithm~\ref{alg:standard_assembly}). In addition, the
\emp{NitscheForm} class gives access to the two overlapping meshes as
well as a quadrature cache to avoid recomputation of quadrature rules.

Degrees of freedom of the cells of $\mesh_0$ that are entirely covered
by the overlapping mesh $\mesh_2$ are \emph{inactive}; that is, those
degrees of freedom do not determine the solution and can be assigned
an arbitrary value. For practical reasons, inactive degrees of freedom
are included in the linear system but are set to zero by inserting
'one' on the diagonal of the corresponding rows and 'zero' in the
right-hand side vector. This is automatically handled by the DOLFIN
function \emp{Matrix::ident\_zeros}.

\section{Numerical examples}
\label{sec:examples}

To demonstrate the efficiency of our realization of the overlapping
mesh method, we study two test examples. The first example is the
Poisson equation. Here, an analytical solution for a suitable source
function allows to verify the implementation by means of convergence
studies. In addition, timings for the computation of both a standard
$P_1$ finite element method and its corresponding Nitsche
approximation are presented. A further breakdown of the timings in the
Nitsche case gives a clear picture of the additional costs associated
with the geometry-related computations and their effect on the overall
computation time. From that perspective, comparing Nitsche with a
simple piecewise linear, continuous finite element method represents
the most challenging test case, since then the extra work required for
geometry-related computations contributes the most to the overall
assembly time.

As a second example, a linear elastic equation is solved with
discontinuous material parameters at the interface between two
overlapping meshes.

\subsection{Poisson equation}

\begin{figure}
  \begin{center}
    \codesize
    \begin{code}
// Function spaces and forms on the overlapped domain
Poisson3D_1::FunctionSpace V1(mesh1);
Poisson3D_1::BilinearForm a1(V1, V1);
Poisson3D_1::LinearForm L1(V1);
L1.f = f; // assign source

// Function spaces and forms on the overlapping domain
Poisson3D_2::FunctionSpace V2(mesh2);
Poisson3D_2::BilinearForm a2(V2, V2);
Poisson3D_2::LinearForm L2(V2);
L2.f = f; // assign source

// Build Nitsche forms
PoissonNitsche::BilinearForm a_nit(a1, a2);
PoissonNitsche::LinearForm b_nit(L1, L2,
                        	 a_nit.overlapping_meshes_ptr(),
                        	 a_nit.overlapped_domain_quadrature_cache_ptr(),
                        	 a_nit.interface_domain_quadrature_cache_ptr());

// Assemble
Matrix A;
Vector b;
NitscheAssembler::assemble(A, a_nit, 0, 0, 0);
NitscheAssembler::assemble(b, b_nit, 0, 0, 0);

// Apply boundary conditions
DirichletBC bc(V1, u0, boundary);
bc.apply(A, b);

// Solve linear system
Vector x;
solve(A, x, b, "cg", "amg_hypre");

// Split solution vector according to domains
Function u1(V1);
Function u2(V2);
a_nit.distribute_solution(x, u1, u2);
    \end{code}
  \end{center}
  \caption{Code example for the Poisson problem. The domain
    decomposition character of Nitsche's method is clearly reflected
    by defining forms on each mesh separately and ``gluing'' them
    together via a~\emp{NitscheForm}. Since both forms are assembled
    into a single matrix, the solution has to be split via the
    \emp{distribute\_solution} function.}
  \label{code:nitsche_problem}
\end{figure}

\begin{figure}
  \begin{center}
    \label{fig:sin_outer_and_inner_mesh}
  \includegraphics[width=0.7\textwidth]{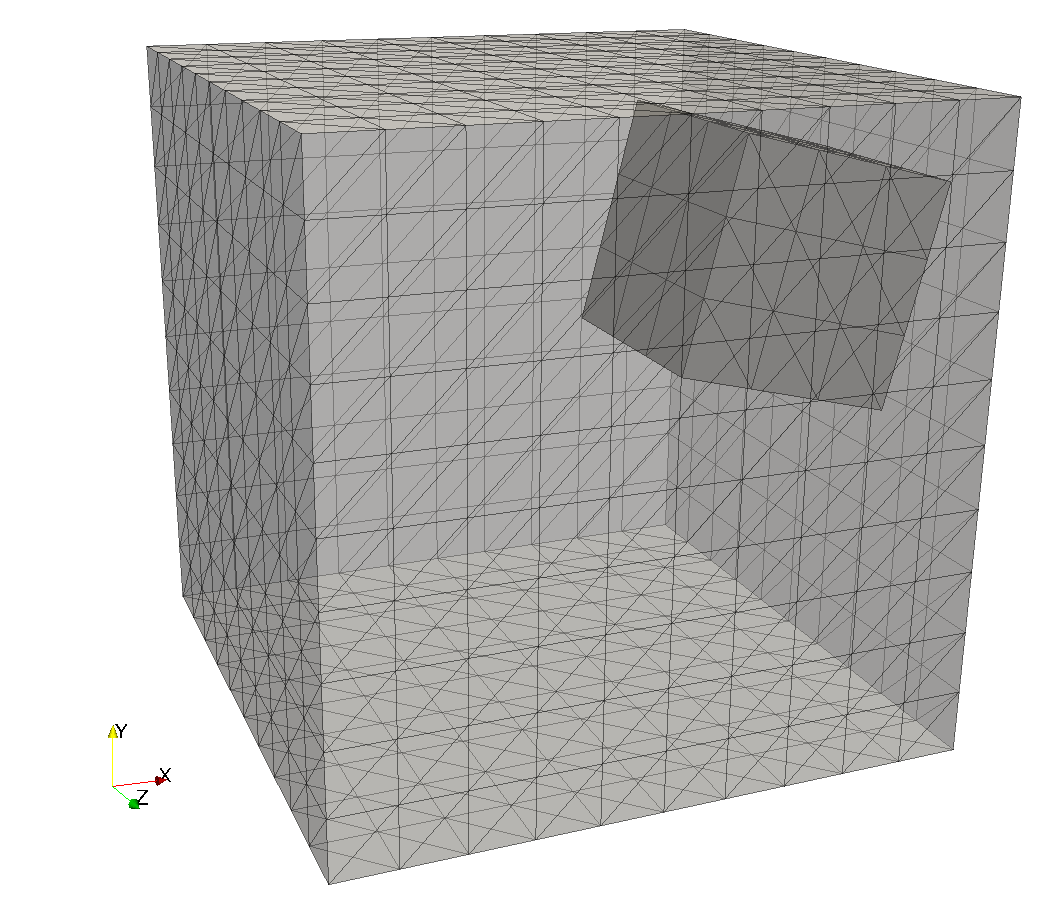}
    \caption{Mesh configuration for the Poisson problem.}
  \end{center}
\end{figure}

We consider the elliptic model problem
\eqref{eq:strong_poisson_equation}--\eqref{eq:neumann_bc} on the
domain $\Omega_0 = (0,1)^3 \subset \R^3$. The overlapping domain
$\Omega_2$ is a translation and rotation of the cube
$\widetilde{\Omega}_2 = [0.3331,0.6669]^3$ according to
Figure~\ref{fig:sin_outer_and_inner_mesh}. The source function $f$
is given by
\begin{equation}
  f(x,y,z) = 3 (2 \pi)^2 \sin(2 \pi x) \sin(2\pi y) \sin(2\pi z).
\end{equation}
Homogeneous Dirichlet boundary conditions are used on the entire
boundary. The exact solution is given by
\begin{equation}
  u(x,y,z) = \sin(2 \pi x) \sin(2\pi y) \sin(2\pi z).
  \label{eq:poisson_exact_solution}
\end{equation}
The penalty parameter $\gamma$ in \eqref{eq:poisson_nitsche_form} is
set to $\gamma = 50$. A code extract from the implementation of the
solver based on DOLFIN-OLM is shown in
Figure~\ref{code:nitsche_problem}.

The Nitsche approximation was computed on a sequence of meshes with
decreasing mesh size $h_{\mathrm{max}} = 1/N$ for the tessellations
$\mesh_0$ and $\mesh_2$, starting at $N = 14$ and stopping at $N =
104$. For the integration over cut cells and interface facets,
barycenter quadrature was employed. The resulting linear systems were
solved by the preconditioned conjugate gradient method as implemented
in PETSc in combination with the algebraic multigrid solver from Hypre
used as a preconditioner. For both the standard finite element method
and the overlapping mesh method, a mesh independent number of CG
iterations was observed (3--4 and 9-11, respectively). All
computations were carried out on a Macbook Pro equipped with a $2.66$
GHz Intel Core i7 processor and $8$ GB of RAM (1066 MHz DDR3). The
benchmarks were repeated 10 times and averaged to obtain the reported
results.

\paragraph{Convergence}

\begin{figure}
  \begin{center}
    \includegraphics[width=0.49\textwidth]{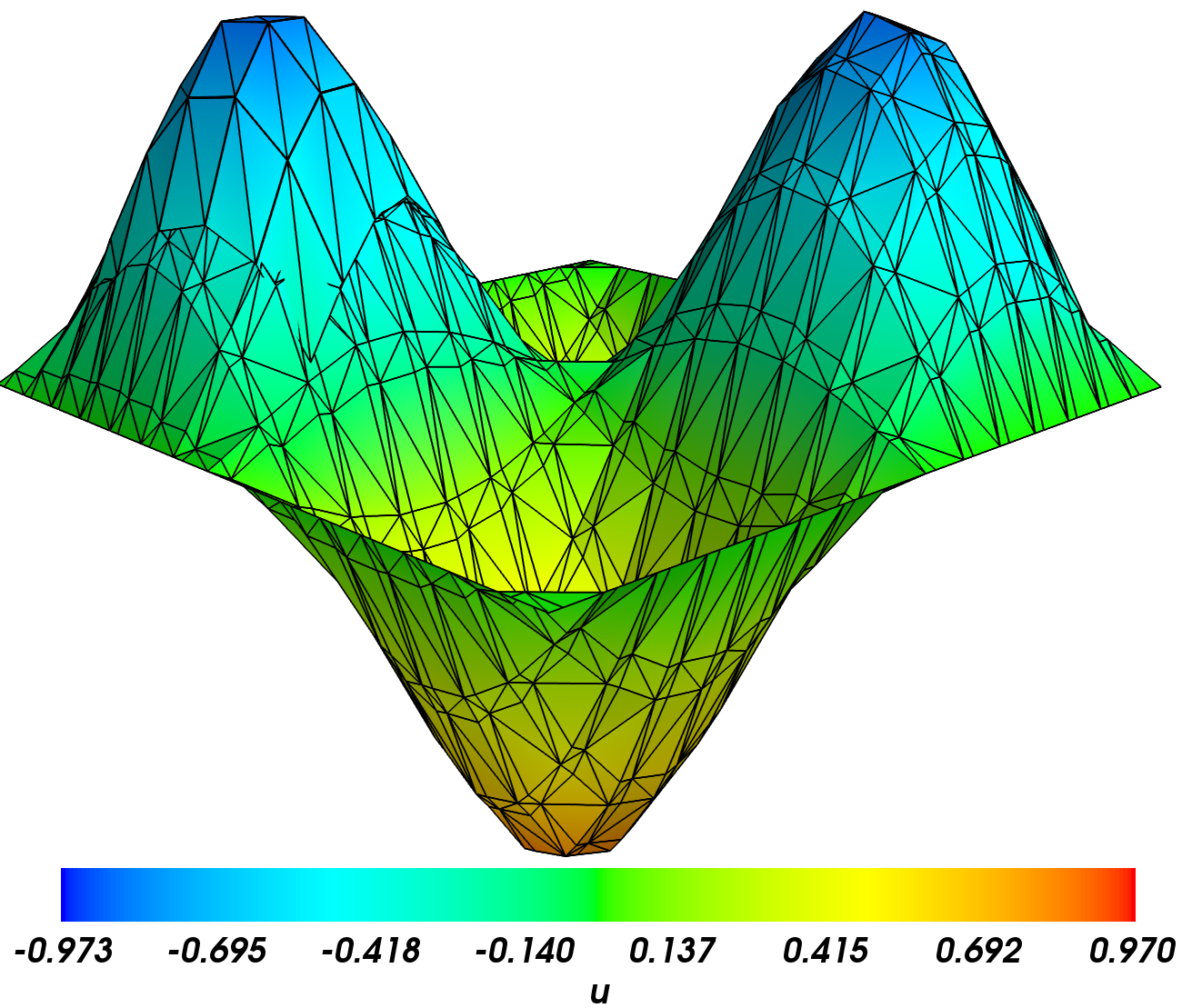}
    \includegraphics[width=0.49\textwidth]{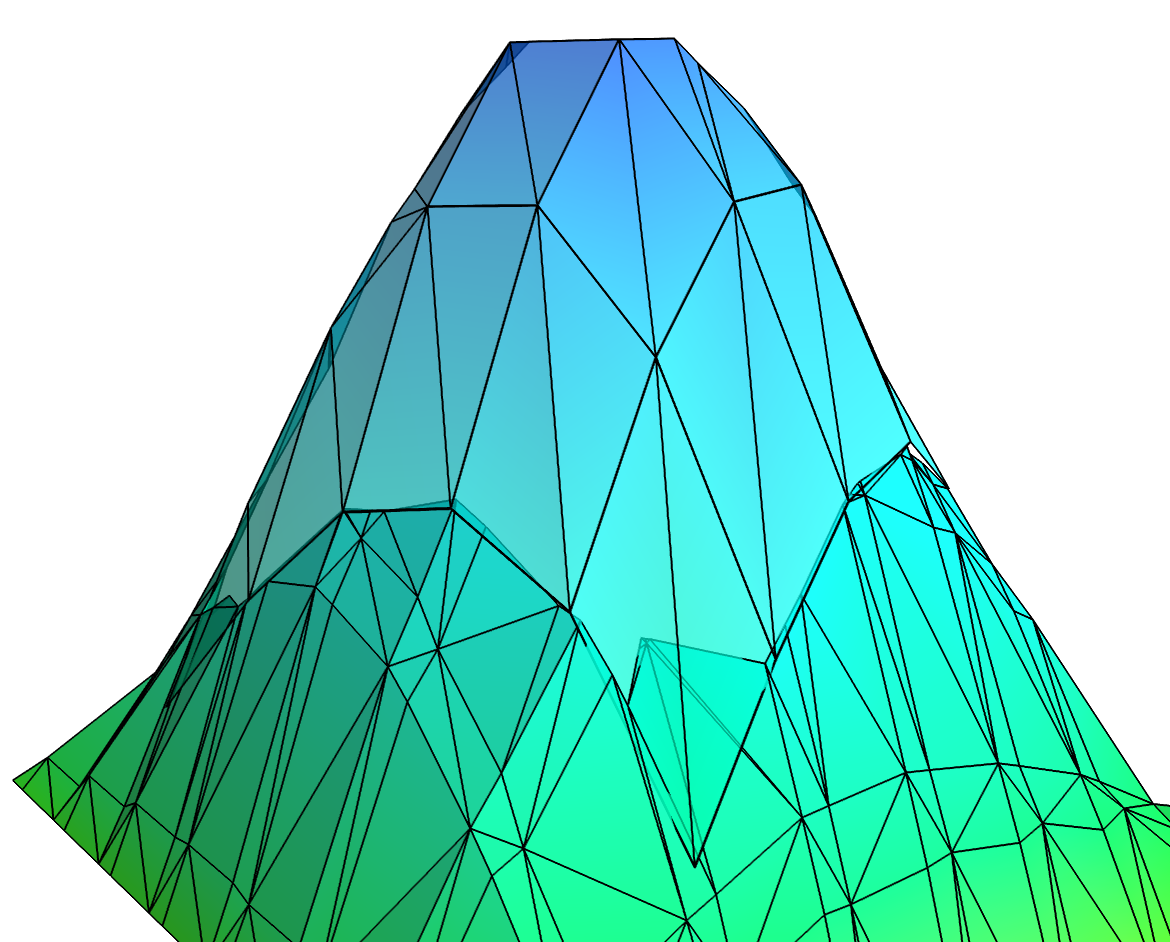}
    \caption{Two-dimensional cross-section of the three-dimensional
      solution of the Poisson problem showing good agreement between
      solutions computed on the overlapped and the overlapping mesh.
      (Left) Both solution parts patched together. The solution on the
      overlapping mesh can be seen far to the left. (Right) A zoom
      revealing the small discontinuity between solutions at the
      common interface.}
    \label{fig:sin_numerical_solution_wo_warp}
  \end{center}
\end{figure}

\begin{figure}
  \begin{center}
    \includegraphics[width=0.95\textwidth]{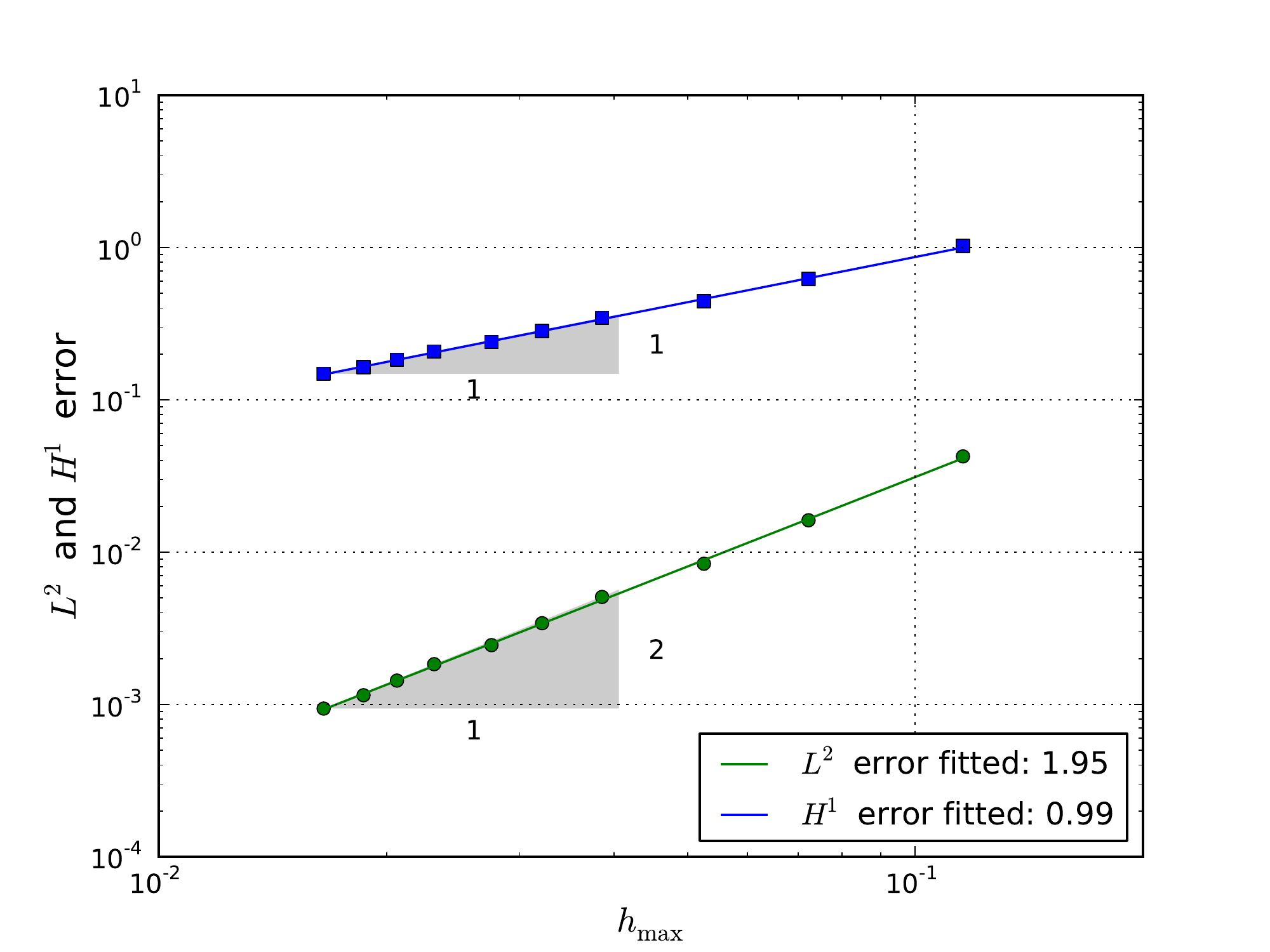}
    \caption{Convergence of the Nitsche approximation in the $L^2$- and
      $H^1$-norm for the Poisson problem.}
      \label{fig:sin_convergence_plot}
  \end{center}
\end{figure}

As the theoretical results recalled in
\eqref{eq:h1norm_error}--\eqref{eq:l2norm_error} predict, an optimal
convergence rate in both in the $H^1$- and $L^2$-norm is observed; see
Figure~\ref{fig:sin_convergence_plot}.
Figure~\ref{fig:sin_numerical_solution_wo_warp} clearly illustrates
the smooth transition from the solution part $u_1$ defined on the
overlapped mesh $\mesh_1$ to the part $u_2$ defined on the overlapping
mesh $\mesh_2$.

\begin{figure}
  \begin{center}
    \includegraphics[height=0.72\textwidth]{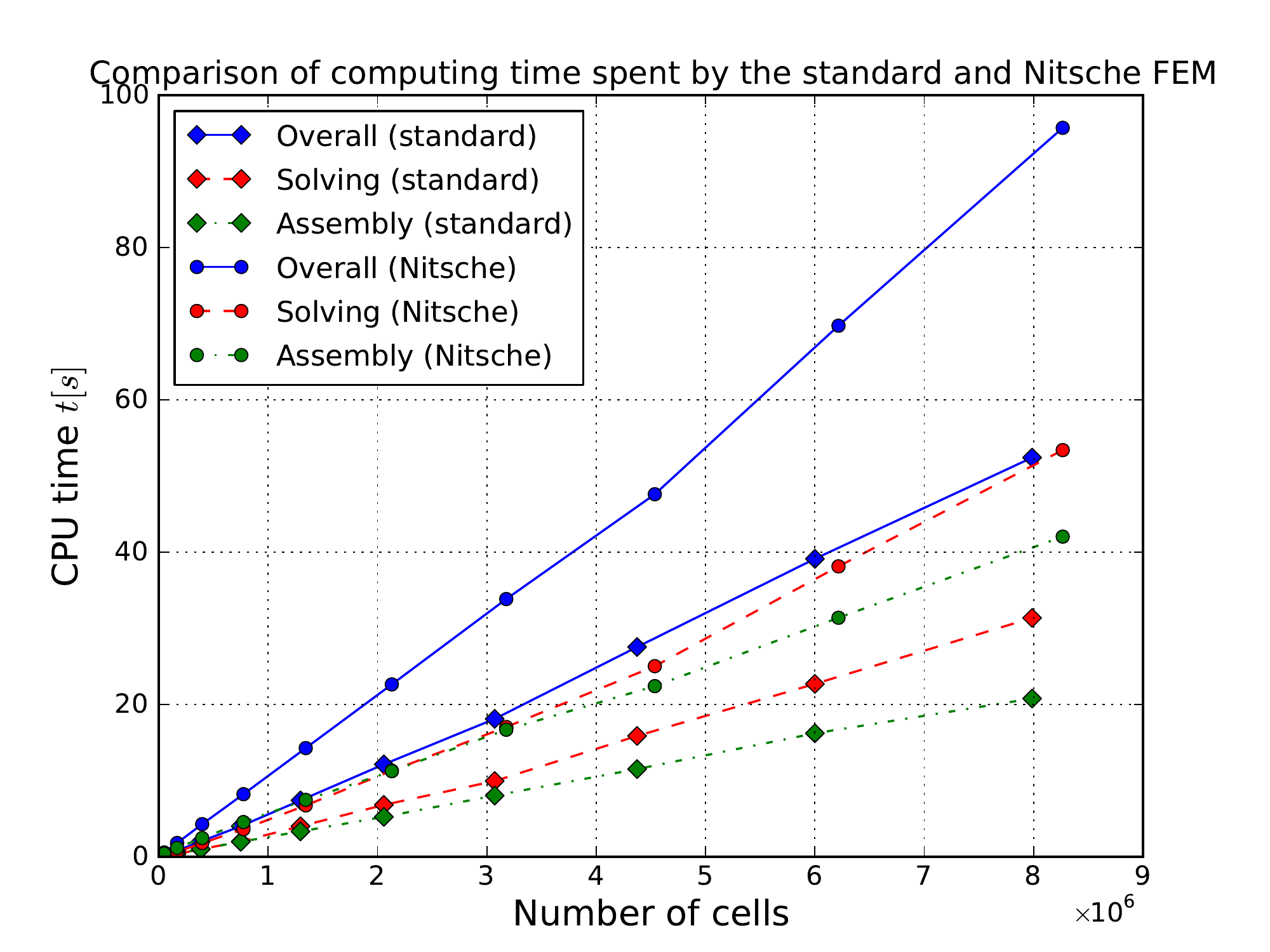}
    \includegraphics[height=0.80\textwidth]{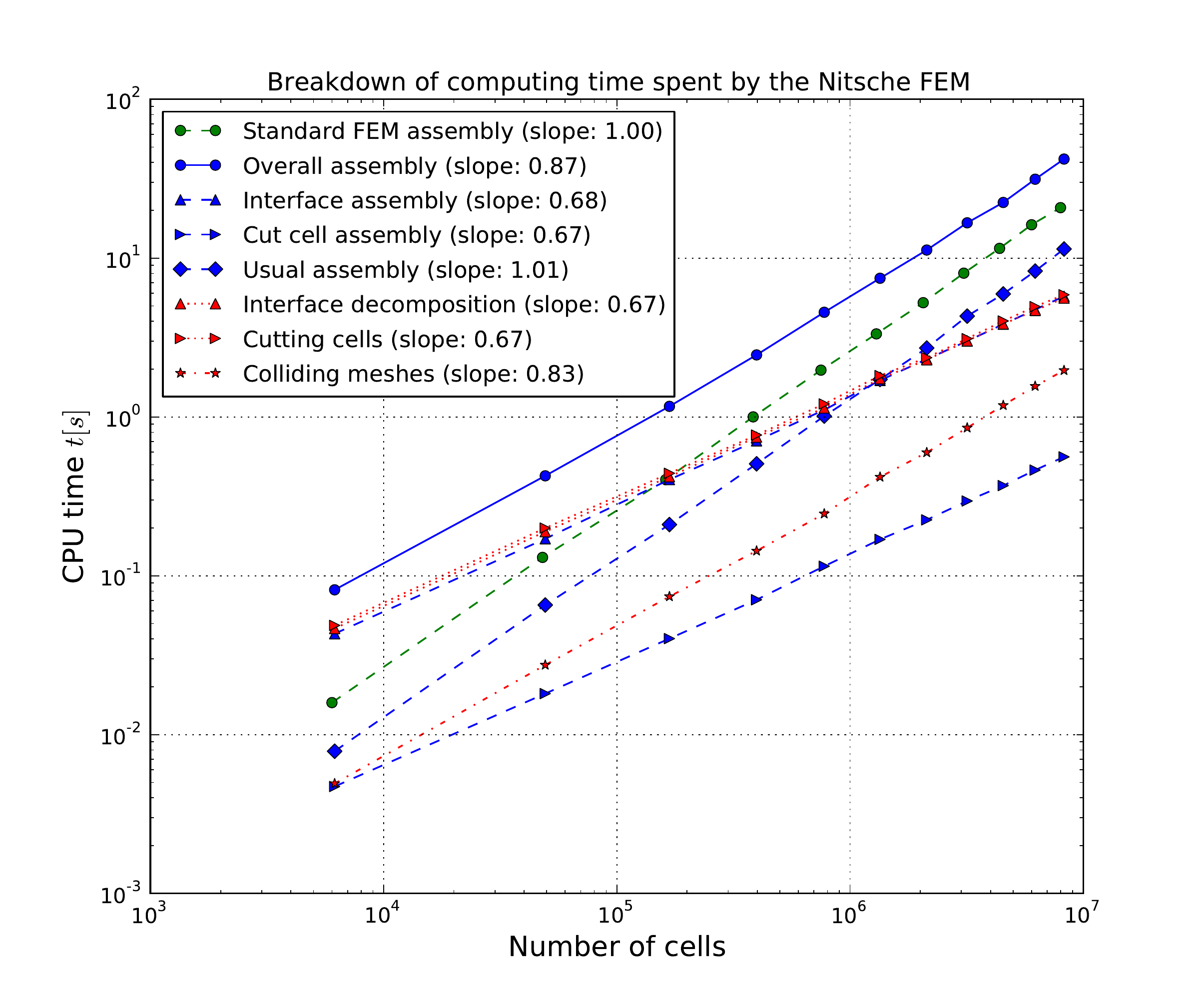}
    \caption{(Top) Overall computing time spent by the standard FEM
      and the Nitsche FEM. (Bottom) Breakdown of the assembly time
      spent by the Nitsche method compared to the overall assembly time
      for the standard FEM. Note the expected slope of 1 for
      the standard (usual) assembly and of $2/3$ for the interface
      related operations. The graph shows that, asymptotically,
      overall assembly time will be dominated by standard assembly.}
    \label{fig:benchmark-overall}
  \end{center}
\end{figure}

\paragraph{Benchmarks}
The comparison in Figure~\ref{fig:benchmark-overall}
shows that the Nitsche overlapping mesh method is only twice as
expensive as the standard finite element method for the same mesh size
and approximately the same number of degrees of freedom. A
breakdown of the computing time shows that both the assembly and the
solution of the linear system become twice as expensive when the
overlapping mesh method is used. The
latter can be attributed to the observed higher iteration numbers. The
CPU time for the linear solve displays a kink in the slope at ca 3 and
4.5 million cells for the standard and Nitsche FEM, respectively,
which may be attributed to the problem size increasing beyond a
hardware-specific threshold (cache or memory size) that induces
overhead.

A further breakdown of the assembly time is shown at the bottom of
Figure~\ref{fig:benchmark-overall}. We note that while the number of
cells scales like $\mathcal{O}(h^{-3})$ on a regular grid, the number
of facets scales like $O(h^{-2})$. Consequently, all purely interface
related operations with linear complexity should scale like $N^{2/3}$,
where $N$ is the number of cells. This is also the case for the
interface decomposition, computation of the cut cells, and integration
over the interface and the cut cells. On the other hand, the initial
collision detection between the overlapped and the overlapping mesh is
somewhat more expensive as it involves the searching of tree-like data
structures (adding a logarithmic factor to the complexity).  However,
this extra cost is negligible compared to the total computing time.
In conclusion, the timings depicted in
Figure~\ref{fig:benchmark-overall} indicate that the cost of assembly
for the Nitsche overlapping mesh method is comparable to standard
finite element assembly, and its relative efficiency increases with
increasing mesh size.

\paragraph{A note on ill-conditioned stiffness matrices}

As analyzed in \citet{Burman2010} and illustrated by numerical
experiments in \cite{Massing2012}, the stiffness matrix
stemming from Nitsche's method on overlapping meshes and related
schemes can be ill-conditioned if the intersection gives very small
elements, compared to the original element size.  Different approaches
have been investigated to cure the schemes from ill-conditioned
systems, either by choosing proper weights for the interface
\citep{JohanssonLarson}, replacing basis functions on small elements
by extensions of basis functions on larger neighboring elements, or
through the use of a so-called ghost-penalty
\citep{Burman2010,BurmanHansbo2011}.
For the current work, the only precaution was to skip all cut cells
with a relative measure of size smaller that $10^{-15}$. In our
numerical experiments, we observed a mesh-independent iteration number
which indicates that ill-conditioning was not present. However, proper
handling of small cells by introducing ghost penalties has been
studied recently for the Stokes problem in \citet{Massing2012b}.

\subsection{Linear elasticity}
As a second example, we consider a linear elastic body occupying a
domain $\Omega_0$ in $\R^3$ consisting of two subdomains $\Omega_1$
and $\Omega_2$ with possibly different material parameters. The
displacement and the stresses are assumed to be continuous across the
interface $\Gamma = \partial \Omega_1 \cap \partial \Omega_2$ between
the subdomains. The corresponding linear elasticity problem then takes
the form: find the displacement $\vect{u}: \Omega_0 \rightarrow \R^3$
such that

\begin{alignat}{3}
  \label{eq:linearelasta}
  -\Div(\bfsigma(\vect{u}_i)) &= \vect{f}_i &\quad &\text{in }  \Omega_i, \quad i = 1,2,
  \\
    \jump{\bfsigma(\vect{u}) \cdot \vect{n}} &= \vect{0} &\quad &\text{on } \Gamma,
    \\
      \jump{\vect{u}} &= \vect{0} &\quad &\text{on } \Gamma,
      \\
      \vect{u} &= \vect{0} &\quad &\text{on } \partial \OmegaD{},
      \\
      \label{eq:linearelastb}
      \bfsigma(\vect{u}) \cdot \vect{n} &= \vect{g} &\quad &\text{on } \partial \OmegaN{}.
\end{alignat}
Here, the stress tensor $\bfsigma$ is related to the displacement
vector $\vect{u}$ by Hooke's law
\begin{equation}
  \bfsigma(\vect{u}_i) = 2 \mu_i \bfepsilon(\vect{u}_i) + \lambda_i
  \text{tr}(\bfepsilon(\vect{u}_i))\mathbf{I} \quad \text{in }
  \Omega_i, \quad i = 1,2,
\end{equation}
where $\lambda_i$ and $\mu_i$ are the Lam\'e parameters in $\Omega_i$
for $i=1,2$ and $\bfepsilon(\vect{v}) = (\Grad(\vect{v}) +
\Grad(\vect{v})^{\top})/{2}$ is the strain tensor.

Nitsche's method for the Poisson equation
\eqref{eq:strong_poisson_equation} proposed
by~\citet{HansboHansboLarson2003} can be adapted to the case of the
linear elastic problem
\eqref{eq:linearelasta}--\eqref{eq:linearelastb} and takes the
following form: find $\vect{u}_h \in \vect{W}_h$ such that
\begin{equation}
  a(\vect{u}_h,\vect{v}) = l(\vect{v}) \quad \foralls \vect{v} \in \vect{W}_h,
\end{equation}
where
\begin{align}
  a(u,v)
  &=
  \sum_{i=1}^2 \int_{\Omega_i} \bfsigma(\vect{u}_i) : \Grad(\vect{v}_i) \dx \\
  &\qquad
  - \underbrace{\int_{\Gamma  }  \bfsigma(\vect{u}_1) \cdot \vect{n} \, \jump{\vect{v}} }_{\text{Stress balance}} \dS
  - \underbrace{\int_{\Gamma}  \bfsigma(\vect{v}_1) \cdot \vect{n} \, \jump{\vect{u}} }_{\text{Symmetrization}} \dS
  + \underbrace{\gamma \int_{\Gamma} h^{-1} \jump{\vect{u}} \cdot \jump{\vect{v}} }_{\text{Penalty/Stabilization}} \dS,
  \\
  l(v) &= \int_\Omega \vect{f}\cdot \vect{v} \dx
         + \int_{\partial \Omega} \vect{g} \cdot \vect{v} \ds,
\end{align}
with $\gamma$ a positive penalty parameter. Assuming that the elastic
material is not nearly incompressible ($\lambda$ remains bounded), we
may extend the analysis in~\citet{HansboHansboLarson2003} and prove
optimal order \apriori{} error estimates. See \citet{HansboLarson2002}
and \citet{BeckerBurmanHansbo2009} for details on discontinuous
Galerkin methods for elasticity problems.

\paragraph{Test configuration and numerical results}

We consider the linear elasticity problem
(\ref{eq:linearelasta})--(\ref{eq:linearelastb}) with $\Omega_0 =
(-2,2)^3\subset \R^3$ which is overlapped by a propeller-like domain
$\Omega_1 = (-1,1) \times P$ where $P = (-1,1) \times (-0.2,0.2) \cup
(-0.2,0.2) \times (-1,1) \subset \R^2$; see
Figure~\ref{fig:intersection_areas}. The right-hand side is zero and
the boundary conditions on $\partial \Omega_0$ are defined by
\begin{equation}
  \begin{cases}
    \vect{u} = 0, &\quad \text{on } (-2,2)^2 \times \{-2\},\\
    \bfsigma(\vect{u}) \cdot \vect{n} = 0, &\quad \text{on } \partial ((-2,2)^2) \times (-2,2),\\
    \bfsigma(\vect{u}) \cdot \vect{n} = \vect{g} &\quad \text{on } (-2,2)^2
    \times \{2\},
  \end{cases}
\end{equation}
where $\vect{g}$ represents a combination of pure tangential,
rotational force and a normal pressure:
\begin{equation}
  \vect{g}(x,y,z) = \dfrac{(-y,x,0)^{\top}}{5\sqrt{x^2 + y^2}} -
  (0,0,2 - \sqrt{x^2 + y^2})^{\top}.
\end{equation}
The Lam\'e parameters are given by $\mu_i = E_i / (2 + 2\nu_i) ,
\lambda_i = E_i \cdot \nu_i /((1 + \nu_i)(1 - 2\nu_i))$ in $\Omega_i$
for $i=1,2,$ with $E_1 = 10$, $E_2 = 0.1\cdot E_1$, $\nu_1 = \nu_2 = 0.3$.

\begin{figure}
  \includegraphics[width=0.49\textwidth]{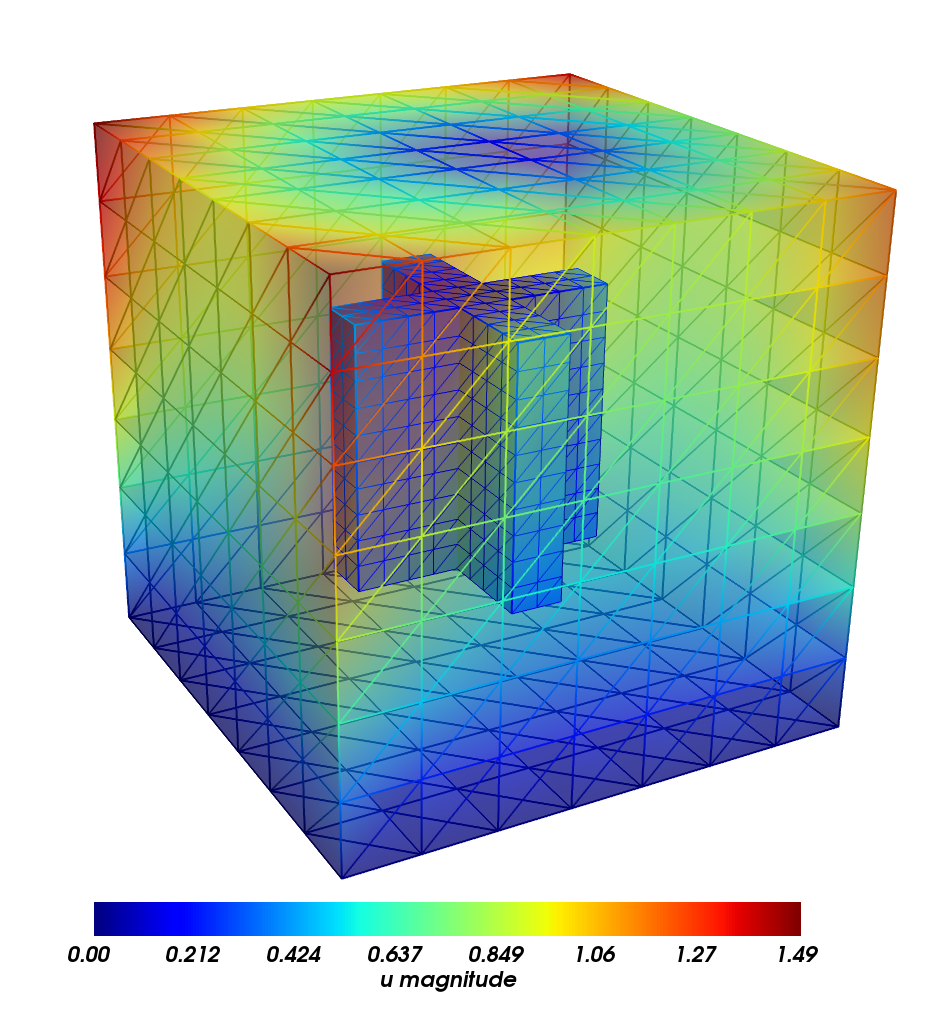}
  \includegraphics[width=0.49\textwidth]{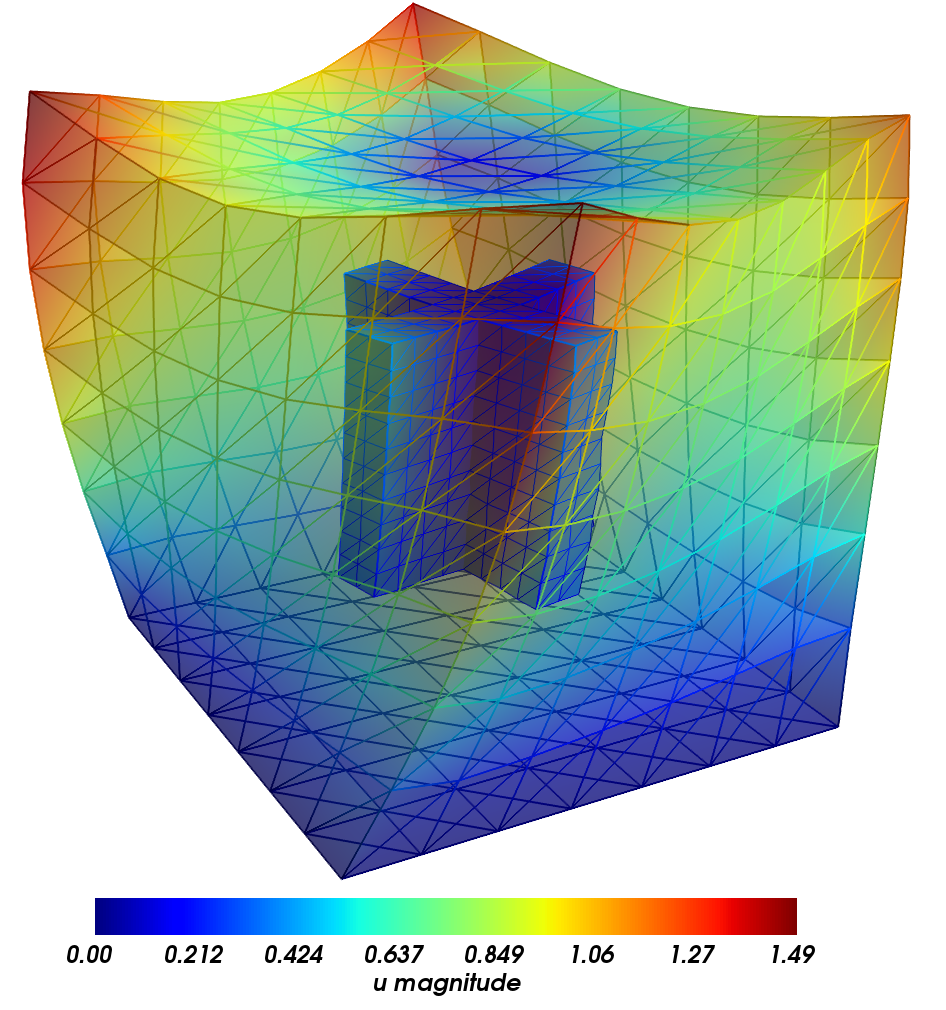}
  \caption{(Left) Undeformed, original domain consisting of a cube
    overlapped by a propeller-like domain. (Right) The deformed
    domain. The color bar corresponds to the norm of the displacement
    $\vect{u}$.}
  \label{fig:soft_material_start_configuration}
\end{figure}

The numerical results are shown in
Figures~\ref{fig:soft_material_start_configuration}
and~\ref{fig:elasticity,magnitude}. The results indicate a ``smooth''
transition of the solution from the overlapped mesh to the overlapping
mesh.

\begin{figure}
  \includegraphics[width=0.48\textwidth]{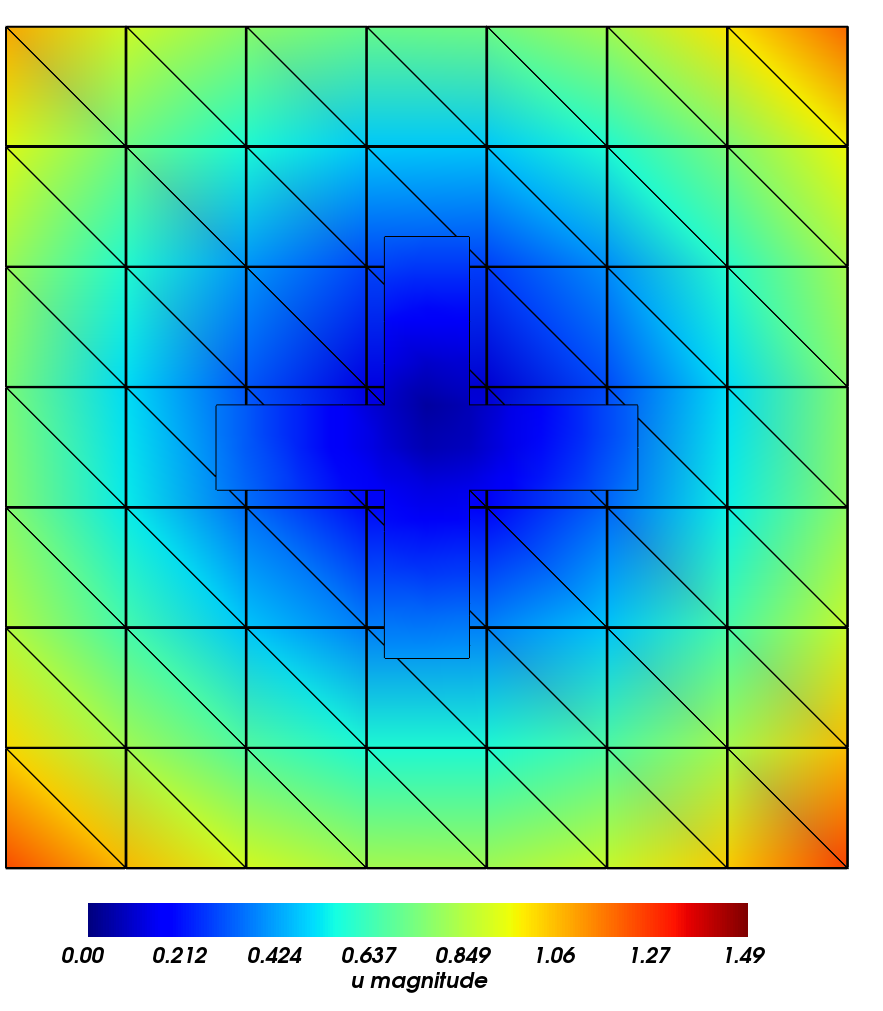}
  \includegraphics[width=0.51\textwidth]{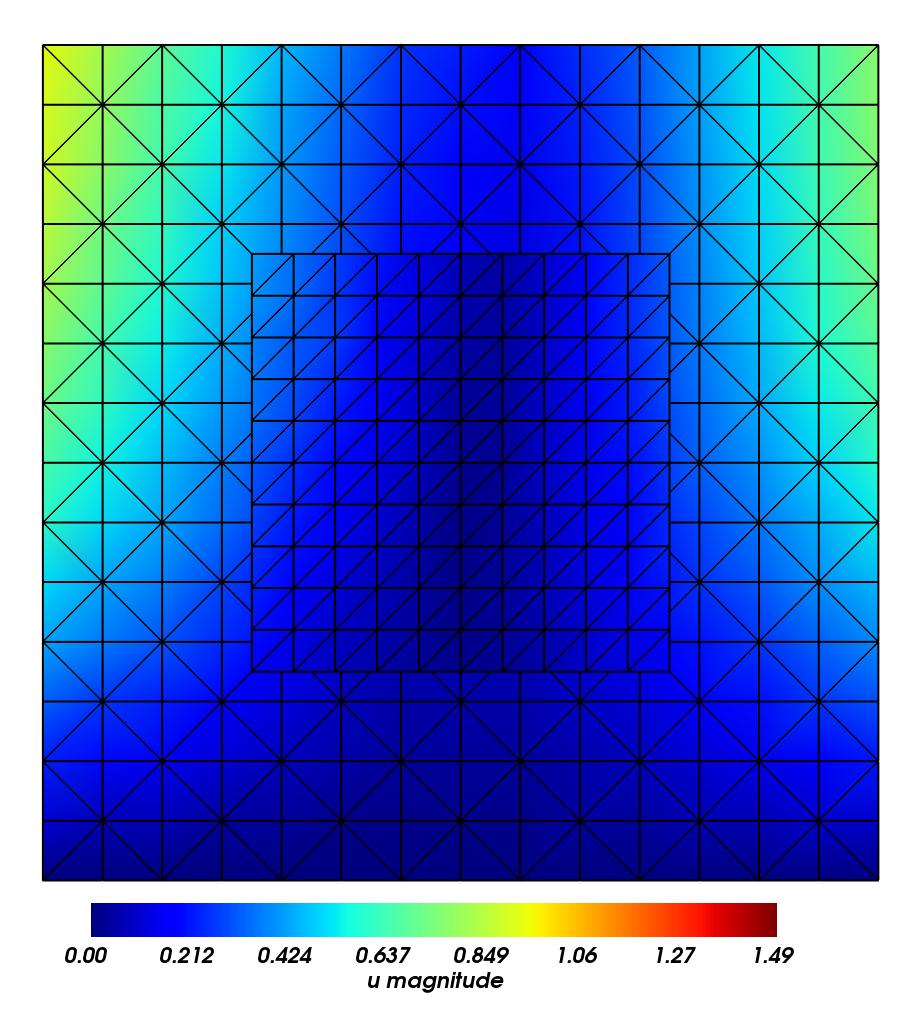}
  \caption{(Left) Magnitude of the displacement $\vect{u}$ in a
    cross-section through the $xy$-plane of the solution from
    Figure~\ref{fig:soft_material_start_configuration}.  (Right)
    Magnitude of the displacement $\vect{u}$ in the
    cross-section through the $xz$-plane. Despite a very coarse resolution, both figures show
    a ``smooth'' transition of the solution from the domain $\Omega_1$
    to $\Omega_2$.}
  \label{fig:elasticity,magnitude}
\end{figure}

\section{Conclusions and outlook}

We have demonstrated that overlapping mesh methods, in particular the
Nitsche overlapping mesh method, may be implemented efficiently in
three space dimensions through the use of tree search data structures
and tools from computational geometry. Numerical tests show that
optimal order convergence is obtained, that the overhead of the
overlapping mesh method compared to a standard finite element method
is small (roughly factor two), and that the overhead is decreasing as
the size of the mesh is increased.

In the future, we plan to extend our implementation and the techniques
studied in this work to handle fluid--structure interaction problems
as well as contact problems. Furthermore, we plan to fully automate
the implementation of Nitsche formulations on overlapping meshes by
adding code generation capabilities for interface terms to FEniCS.

\section*{Acknowledgments}

We thank Harish Narayanan for valuable discussions regarding the
elasticity test case. This work is supported by an Outstanding Young
Investigator grant from the Research Council of Norway, NFR
180450. This work is also supported by a Center of Excellence grant
from the Research Council of Norway to the Center for Biomedical
Computing at Simula Research Laboratory.

\bibliographystyle{plainnat}
\bibliography{bibliography}

\end{document}

%% file: pdf/intersected_cell_with_dofs_1.pdf_tex
\begingroup%
  \makeatletter%
  \providecommand\color[2][]{%
    \errmessage{(Inkscape) Color is used for the text in Inkscape, but the package 'color.sty' is not loaded}%
    \renewcommand\color[2][]{}%
  }%
  \providecommand\transparent[1]{%
    \errmessage{(Inkscape) Transparency is used (non-zero) for the text in Inkscape, but the package 'transparent.sty' is not loaded}%
    \renewcommand\transparent[1]{}%
  }%
  \providecommand\rotatebox[2]{#2}%
  \ifx\svgwidth\undefined%
    \setlength{\unitlength}{510.11209611bp}%
    \ifx\svgscale\undefined%
      \relax%
    \else%
      \setlength{\unitlength}{\unitlength * \real{\svgscale}}%
    \fi%
  \else%
    \setlength{\unitlength}{\svgwidth}%
  \fi%
  \global\let\svgwidth\undefined%
  \global\let\svgscale\undefined%
  \makeatother%
  \begin{picture}(1,0.64143658)%
    \put(0,0){\includegraphics[width=\unitlength]{intersected_cell_with_dofs_1.pdf}}%
    \put(0.48259623,0.43964868){\color[rgb]{0,0,0}\makebox(0,0)[lb]{\smash{$T^2_k$}}}%
    \put(0.47534723,0.29005741){\color[rgb]{0,0,0}\makebox(0,0)[lb]{\smash{$\Gamma_{kl}$}}}%
    \put(0.59577065,0.29008926){\color[rgb]{0,0,0}\makebox(0,0)[lb]{\smash{$\Gamma_{kl'}$}}}%
    \put(0.40332233,0.11029331){\color[rgb]{0,0,0}\makebox(0,0)[lb]{\smash{$P^0_l = T^0_l \cap \Omega_1$}}}%
    \put(0.72481391,0.15618162){\color[rgb]{0,0,0}\rotatebox{60.22022482}{\makebox(0,0)[lb]{\smash{$P^0_{l'} = T^0_{l'} \cap \Omega_1$}}}}%
    \put(0.64453602,0.44014288){\color[rgb]{0,0,0}\makebox(0,0)[lb]{\smash{$T^2_{k'}$}}}%
    \put(0.69893183,0.31444594){\color[rgb]{0,0,0}\rotatebox{52.61248982}{\makebox(0,0)[lb]{\smash{$\Gamma_{k'l'}$}}}}%
  \end{picture}%
\endgroup%